\documentclass[oneside, a4paper, 12pt]{amsart}
\pdfoutput=1

\usepackage{xparse}
\usepackage{environ}
\usepackage[a4paper, margin=3cm]{geometry}
\usepackage{adjustbox}
\usepackage[utf8]{inputenc}
\usepackage[T1]{fontenc}
\usepackage[sb]{libertinus}
\usepackage[amsthm,lcgreekalpha]{libertinust1math}
\usepackage[cal=boondoxo, calscaled=0.96]{mathalfa}
\usepackage{textcomp}
\usepackage[tracking=true,kerning=true]{microtype}
\SetTracking% Slightly space out small caps
{encoding=T1, shape=sc}
{10}
\SetProtrusion% Protrude sub- and superscripts at the right margin
{encoding=T1,size={7,8}}
{1={ ,750},2={ ,500},3={ ,500},4={ ,500},5={ ,500},
  6={ ,500},7={ ,600},8={ ,500},9={ ,500},0={ ,500}}
\usepackage{mathtools}
\mathtoolsset{mathic=true}

\usepackage[english]{babel}
\usepackage{csquotes}
\usepackage{CJKutf8}

\usepackage{hyphenat}
\usepackage{tikz}
\usepackage{amssymb}
\usetikzlibrary{backgrounds, calc, cd}

\pgfdeclarelayer{edgelayer}
\pgfdeclarelayer{nodelayer}
\pgfsetlayers{background,edgelayer,nodelayer,main}
\tikzset{baseline=(current bounding box.center)}
\tikzset{label distance=-0.15cm}
\tikzset{font=\scriptsize}
\tikzstyle{none} = [inner sep=3pt,outer sep=0pt]

\tikzset{curve/.style={settings={#1},to path={(\tikztostart)
      .. controls ($(\tikztostart)!\pv{pos}!(\tikztotarget)!\pv{height}!270:(\tikztotarget)$)
      and ($(\tikztostart)!1-\pv{pos}!(\tikztotarget)!\pv{height}!270:(\tikztotarget)$)
      .. (\tikztotarget)\tikztonodes}},
  settings/.code={\tikzset{quiver/.cd,#1}
    \def\pv##1{\pgfkeysvalueof{/tikz/quiver/##1}}},
  quiver/.cd,pos/.initial=0.35,height/.initial=0}
\tikzset{tail reversed/.code={\pgfsetarrowsstart{tikzcd to}}}
\tikzset{2tail/.code={\pgfsetarrowsstart{Implies[reversed]}}}
\tikzset{2tail reversed/.code={\pgfsetarrowsstart{Implies}}}
\tikzset{no body/.style={/tikz/dash pattern=on 0 off 1mm}}
\tikzcdset{
  arrow style=tikz,
  diagrams={>={Straight Barb}}
}

\definecolor{red}{RGB}{175, 49, 39}             %#AF3127
\usepackage{xspace}
\usepackage[normalem]{ulem}
\usepackage{enumitem}
\usepackage{aligned-overset}
\usepackage[square,numbers,sort&compress]{natbib} % Bibliography
\usepackage[colorlinks, linktoc=page, linkcolor={red}, citecolor={red}, urlcolor={red}]{hyperref}

\newlist{thmlist}{enumerate}{1}
\setlist[thmlist]{label=(\roman{thmlisti}), ref=\thetheorem.(\roman{thmlisti}), noitemsep}
\newlist{thmlist*}{enumerate}{1}
\setlist[thmlist*]{label=(\roman{thmlist*i}), ref=(\roman{thmlist*i}), noitemsep}

\numberwithin{equation}{section}

\DeclareRobustCommand{\SkipTocEntry}[5]{}

\theoremstyle{plain}
\newtheorem{theorem}{Theorem}
\newtheorem{lemma}[theorem]{Lemma}
\newtheorem{corollary}[theorem]{Corollary}
\newtheorem{proposition}[theorem]{Proposition}
\newtheorem{hypothesis}[theorem]{Hypothesis}

\theoremstyle{definition}
\newtheorem{definition}[theorem]{Definition}
\newtheorem{example}[theorem]{Example}
\newtheorem{remark}[theorem]{Remark}

\newcommand*{\eg}{e.g.\xspace}
\newcommand*{\ie}{i.e.\xspace}

\DeclareRobustCommand{\neghairskip}{\hskip -0.075em\relax}
\newcommand*{\GV}{Grothendieck--\neghairskip{}Verdier\xspace}

\renewcommand*{\colon}{\!\nobreak\mskip2mu\mathpunct{}\nonscript%
  \mkern-\thinmuskip{:}\mskip6muplus1mu\relax%
}
\newcommand*{\from}{\ensuremath{\colon}} % As in \(f\from A \to B\).

\newcommand*{\yo}{\text{\begin{CJK}{UTF8}{min}よ\end{CJK}}} % Yoneda embedding
\renewcommand*{\to}{\ensuremath{\longrightarrow}}
\renewcommand*{\mapsto}{\ensuremath{\longmapsto}}
\newcommand*{\iso}{\overset{\smash{\raisebox{-0.65ex}{\ensuremath{\scriptstyle\sim}}}\,}{\to}}

\newcommand*{\cat}[1]{\ensuremath{\mathcal{#1}}} % A category
\newcommand*{\ConcreteCat}[1]{\ensuremath{\mathsf{#1}}\xspace}
\newcommand*{\mky}{\ConcreteCat{Mky}}
\newcommand*{\mkyfin}{\ConcreteCat{mky}}
\newcommand*{\kvect}[1][\Bbbk]{\ensuremath{\mathsf{vect}_{#1}}}
\newcommand*{\kVect}[1][\Bbbk]{\ensuremath{\mathsf{Vect}_{#1}}}
\newcommand*{\Sp}[1]{\ensuremath{\mathbb{S}\mathsf{p}_{#1}}\xspace}
\newcommand*{\Cpsh}[1][C]{\ensuremath{\widehat{\cat{#1}}}}
\newcommand{\slcrys}{\ensuremath{\mathfrak{sl}_{2}(\mathbb{C})\textnormal{-}\mathbf{crys}}}
\DeclareMathOperator*{\Ob}{Ob}

\makeatletter
\newcommand{\lmod}[1]{#1\textnormal{-mod}}

\newcommand{\rmod}[1]{\textnormal{mod-}#1}
\newcommand{\rMod}{\@ifstar{\ST@rMd}{\ST@@rMd}}
\newcommand\ST@rMd[1]{\ensuremath{{\cat{M}}_{#1}}}
\newcommand\ST@@rMd[1]{\textnormal{Mod-}#1}
\makeatother

\newcommand{\rcomod}[1]{\textnormal{comod-}#1}

\makeatletter
\newcommand*{\@smallTriangle}[2]{\vcenter{\hbox{\scalebox{0.75}{\ensuremath{#1#2}}}}}
\newcommand*{\@smallerTriangle}[2]{\vcenter{\hbox{\scalebox{0.60}{\ensuremath{#1#2}}}}}
\newcommand*{\lact}{\mathbin{\mathpalette\@smallTriangle\rhd}}
\newcommand*{\ract}{\mathbin{\mathpalette\@smallTriangle\lhd}}
\newcommand*{\blact}{\mathbin{\mathpalette\@smallerTriangle\blacktriangleright}}
\makeatother

\DeclareMathOperator*{\soc}{soc}
\newcommand*{\id}{\mathrm{id}}
\newcommand*{\Id}{\mathrm{Id}}

\newcommand*{\cocolon}{\nobreak\mskip6mu plus1mu\mathpunct{}\nonscript\mkern-\thinmuskip{:\!}\mskip2mu\relax}
\newcommand*{\adj}[4]{\ensuremath{#1 \colon #3 \rightleftarrows #4 \cocolon #2}\xspace}
\newcommand*{\adjspace}{\,}
\newcommand*{\adjoint}{\adjspace\dashv\adjspace}  % L ⊣ R
\newcommand*{\lad}{\adjoint}

\newcommand*{\ld}[1]{\prescript{\vee}{}{\!#1}}
\newcommand*{\ev}{\mathrm{ev}}
\newcommand*{\coev}{\mathrm{coev}}
\DeclareMathOperator*{\tr}{tr}

\DeclareMathOperator{\im}{im}        % The image
\DeclareMathOperator{\Hom}{Hom}      % The homomorphisms
\DeclareMathOperator{\End}{End}      % The endomorphisms
\DeclareMathOperator{\Aut}{Aut}      % The automorphisms
\DeclareMathOperator{\Rep}{Rep}
\DeclareMathOperator{\GL}{GL}
\DeclareMathOperator{\PGL}{PGL}
\DeclareMathOperator{\spanset}{span}
\DeclareMathOperator{\op}{op}     % The opposite
\DeclareMathOperator{\characteristic}{char} % The characteristic
\DeclareMathOperator{\fin}{fin}             % Finite-dimensional
\DeclareMathOperator{\can}{can}             % A canonical morphism
\DeclareMathOperator{\colim}{colim}         % The colimit
\newcommand*{\blank}{{-}}                   % A unspecified variable
\newcommand*{\bblank}{{=}}
\newcommand*{\tensorop}{\otimes^{\op}}               % An alias for the opposite tp
\renewcommand*{\k}{\ensuremath{\mathbb{k}}\xspace}

\makeatletter
\def\rd{\@ifstar\@rdstar\@rdnostar}
\def\@rdstar#1{\ensuremath{{#1}^{\ast}}}
\def\@rdnostar#1{\ensuremath{{#1}^{\vee}}}
\def\rrd{\@ifstar\@rrdstar\@rrdnostar}
\def\@rrdstar#1{\ensuremath{{#1}^{\ast\ast}}}
\def\@rrdnostar#1{\ensuremath{{#1}^{\vee\vee}}}

\newcommand{\ostar}{\mathbin{\mathpalette\make@circled\star}}
\newcommand{\make@circled}[2]{%
  \ooalign{\(\m@th#1\smallbigcirc{#1}\)\cr\hidewidth\(\m@th#1#2\)\hidewidth\cr}%
}
\newcommand{\smallbigcirc}[1]{%
  \vcenter{\hbox{\scalebox{1.3}{\(\m@th#1\bigcirc\)}}}%
}
\makeatother

\DeclareFontFamily{U}{DSSerif}{\skewchar \font =45}% openface
\DeclareFontShape{U}{DSSerif}{m}{n}{<-> s*[1]  DSSerif}{}
\DeclareFontSubstitution{U}{DSSerif}{m}{n}
\DeclareMathAlphabet{\mathbbbb}{U}{DSSerif}{m}{n}
\newcommand*{\1}{\text{\usefont{U}{DSSerif}{m}{n}1}}

\makeatletter
\newtheorem*{rep@theorem}{\rep@title}
\newenvironment{rtheorem}[1]{\def\rep@title{Theorem~\ref{#1}}\begin{rep@theorem}}{\end{rep@theorem}}
\makeatother

\author{Sebastian Halbig}
\address{S.H., Philipps-Universität Marburg, Arbeitsgruppe Algebraische Lie-Theorie, Hans-Meer\-wein-Straße 6, 35043 Marburg}
\email{sebastian.halbig@uni-marburg.de}

\author{Tony Zorman}
\address{T.Z., Universität Hamburg, Fachbereich Mathematik, Bereich Algebra und Zahlentheorie, Bundesstraße 55, D-20146 Hamburg, Germany}
\email{tony.zorman@uni-hamburg.de}

\date{\today}
\subjclass[2020]{18D15(primary), 18M10(secondary)}
\keywords{closed monoidal categories, rigid monoidal categories, autonomous categories, \GV{} categories, \(*\)-autonomous categories}
\title{Duality in Monoidal Categories}
\begin{document}

\begin{abstract}
  We compare closed and rigid monoidal categories.
  Closedness is defined by the tensor product having a right adjoint: the internal hom functor.
  Rigidity, on the other hand, generalises the duality of finite-dimensional vector spaces.
  In the latter, the internal hom functor is implemented by tensoring with the respective duals.
  This raises the question: can one decide whether a closed monoidal category is rigid, simply by verifying that the internal hom is tensor-representable?
  We provide a counterexample in terms of the category of \(\mathfrak{sl}_2(\mathbb{C})\)-crystals.
  As a byproduct, we obtain characterisations of the \GV{} duality and rigidity of functor categories endowed with Day convolution as their tensor product.
  This has various applications, three of which we study in detail:
  generalisations of quasi-Frobenius algebras, called QF-2 algebras;
  Mackey functors, where we prove that, as expected due to work of Bouc, an object being rigidly dualisable is equivalent to it being finitely-generated projective;
  and crossed modules of finite groups, where we associate to each of these objects a \GV{} category of group-graded representations.
\end{abstract}

\maketitle

% \microtypesetup{protrusion=false}
\tableofcontents
% \microtypesetup{protrusion=true}

\section{Introduction}\label{sec:intro}

Monoidal categories are a ubiquitous tool in mathematics, physics, and computer science.
They often come equipped with additional structures, such as braidings or twists~\cite{Baez2011}.
The aim of this paper is to compare several notions of \emph{monoidal dualities}.
On one hand, we have closedness of monoidal categories, and on the other hand lies rigidity,
generalising the duals of finite-dimensional vector spaces.
\GV{} categories~\cite{boyarchenko13:groth-verdier}, which are also called \(*\)-autonomous categories~\cite{barr79}, describe a duality theory between the strict confinements of rigidity,
and the very general notion of monoidal closedness.
Roughly speaking, these are pairs comprising a monoidal category and a chosen object, called the dualising object, such that a functor related to the internal hom is invertible.
While closed and rigid categories have left and right variants, \GV{} duality is an inherently ambidextrous notion.

One goal of the present article is to study whether one can decide
if a closed monoidal category is rigid
by verifying that the internal hom is \emph{tensor-representable}: given by tensoring with an object.
The authors were made aware of this question by Chris Heunen during the
\href{https://conferences.leeds.ac.uk/bcqt2022/}{\textsc{bcqt}\oldstylenums{2022}} summer school,
who suggested that this is not the case in general.

Section~\ref{sec:dualities} provides a brief summary of the relevant notions concerning rigidity and closed monoidal structures.

In Section~\ref{sec:tensor-rep->gv}, we develop the notion of tensor-representable categories; \ie, categories whose internal hom is given by tensoring with `duals'.
This comes in two versions—tensoring from the left or from the right—and we show that if certain natural transformations relating these two variants are invertible, tensor-representability entails a special form of  \GV{} duality with the monoidal unit as dualising object, see Propositions~\ref{prop:trep-criteria} and~\ref{prop:right-tensor-rep->right-gv}.

In Section~\ref{sec:tensor-rep-functor-cats} we investigate closed and \GV{} structures on (subcategories of) functor categories endowed with the Day convolution as its tensor product.
In particular, in Corollary~\ref{cor:pres-ref-dual-struct} we establish that the Cauchy completion \(\overline{ \cat{C}}\) of (the opposite of) a \(\Bbbk\)-linear closed category \(\cat{C}\) inherits the same type of duality as \( \cat{C}^{\op}\).

We apply this abstract machinery  to representation-theoretic examples in Section~\ref{sec:mky-2-grps}.
First, we consider a variant of the Temperley--Lieb category \(\mathsf{TL}\) whose Cauchy completion can be identified with the category of \(\mathfrak{sl}_{2}(\mathbb{C})\)-crystals, \cite{alqady-stroiński2025:TemperleyLieb}.
It is a non-rigid \GV{} category, see \cite{etingof-penneys2025:RigidityNonNeglegibleModerateGrowth,alqady-stroiński2025:TemperleyLieb}, and we show:
\begin{rtheorem}{thm:sl2-is-trep}
  The category \(\slcrys\) of \(\mathfrak{sl}_2(\mathbb{C})\)-crystals is a non-rigid tensor\hyp{}representable Grothendieck--Verdier category and its dualising functor satisfies \(D_{r}^{2}\cong \Id_{r}\).
\end{rtheorem}
As a result, we obtain a strict hierarchy
\[
  \text{Rigid} \subsetneq \text{Tensor-representable}\footnotemark \subsetneq \text{\GV} \subsetneq \text{Closed}.
\]
\footnotetext{\,%
  Here, we assume that the natural transformations `comparing' the left and right-handed version of tensor-representability are invertible;
  see in particular Propositions~\ref{prop:trep-criteria} and~\ref{prop:right-tensor-rep->right-gv}.
  For rigid categories this is always the case.%
}%
For the various inclusions, see—from left to right—%
Theorems~\ref{thm:compact-closed->closed} and~\ref{thm:sl2-is-trep},
a combination of Propositions~\ref{prop:trep-criteria} and \ref{prop:right-tensor-rep->right-gv} and conversely Example~\ref{ex: r-cat-not-always-trep},
as well as Proposition~\ref{prop:gv->closed} and Remark~\ref{rmk:closed-but-not-gv}.

We furthermore investigate three examples of abelian closed monoidal functor categories.
First, finite Boolean algebras.
These are, roughly speaking, finite sets endowed with the operations `and', `or', and `negation'.
Functors from the associated poset-category to finite-dimensional vector spaces can be identified with modules over an associative and unital algebra.
Negation induces a \GV{} duality on this category and allows us to find a correspondence between projective and injective modules.
This is discussed in Proposition~\ref{prop:qf2-algebras-from-boolean} where we furthermore show that in general these algebras are not quasi-Frobenius—their projectives need not be injective—but satisfy the QF-2 property.
That is, every irreducible projective has a simple socle.
Next, we study finite-dimensional Mackey functors;
in~\cite[Lemma~2.2]{bouc05:mackey}, Bouc sketched an argument that the rigidly dualisable objects in this category are precisely the finitely-generated projective objects.
We prove this result with an alternative method in Proposition~\ref{prop:mky-GV-and-rigid}.
The article is concluded by discussing finite strict 2-groups.
They arise from crossed modules and induce a category of group-graded representations of a finite group.
We show in Propositions~\ref{prop:2-groups-tensor-product}~and~\ref{prop:2-grp-GV-and-rigid} that it is an \(r\)-category,
whose rigidity depends on the characteristic of the ground field.

\addtocontents{toc}{\SkipTocEntry}
\subsection*{Acknowledgements}

We would like to thank Chris Heunen, JS Lemay, Ross Street, and Mateusz Stroiński for fruitful discussions,
and Joseph Martin for pointing out several typos in the first version of this document.
In particular, we want to thank Pavel Etingof for raising the question whether \(\mathfrak{sl}_{2}\)-crystals can be tensor-representable.
Furthermore, we thank the anonymous referees for their invaluable feedback regarding the presentation, organisation, and clarity of the work.
T.Z.~was supported by the DFG grant KR 5036/2--1 and acknowledges support from the Deutsche Forschungsgemeinschaft -- SFB 1624 -- ``Higher structures, moduli spaces and integrability'' -- 506632645.

\numberwithin{theorem}{section}
\section{A ladder of duality concepts for monoidal categories}\label{sec:dualities}
We assume the reader's familiarity with standard concepts of category theory; in particular, adjunctions and monoidal categories as discussed for example in~\cite{MacLane1998} and~\cite{Etingof2015}.
Since the duality concepts we are going to consider are preserved, as well as reflected, by monoidal equivalences—see~\cite{lindner78:adjun}—we restrict ourselves to the strict setting.
As such, let \(\cat{C}\) be a strict monoidal category with \(\blank \otimes \blank \from \cat{C} \times \cat{C} \to \cat{C}\) as its \emph{tensor product} and \(1 \in \cat{C}\) as its \emph{unit}.
\bigskip

Let \(x \in\cat{C}\). Its  \emph{right internal hom}, if it exists, is a functor \({[x, \blank]}_r \from \cat{C} \to \cat{C}\), such that
\begin{equation}\label{eq:right-internal hom-def}
  \adj{x \otimes \blank}{{[x, \blank]}_r}{\cat{C}}{\cat{C}}.
\end{equation}
This adjunction is referred to as the (right) \emph{tensor-hom adjunction}.
In case, the right internal hom exists for all objects, we call \(\cat{C}\) a \emph{right closed category}.
The object-wise tensor-hom adjunctions of a closed category specify a unique functor \({[\blank, \blank]}_r \from \cat{C}^{\op}\times\cat{C} \to\cat{C}\) such that we have \(\cat{C}(x \otimes y, z) \cong \cat{C}(y, {[x,z]}_r)\) natural in all three variables, see~\cite[Section IV.7]{MacLane1998}.

We write \(\cat{C}^{\tensorop}\) for the category \(\cat{C}\) equipped with its \emph{opposite tensor product} \(\tensorop\), given by \(x \tensorop y \eqdef y \otimes x\) on objects.
Analogous considerations to the previous paragraph for tensoring with a fixed object on the right lead to the notion of \emph{left closedness},
which can be succinctly characterised by \(\cat{C}^{\tensorop}\) being right closed.
Provided \(\cat{C}\) is both left and right closed, we will simply refer to it as a \emph{closed} or \emph{biclosed} monoidal category,
where the latter may be used to emphasise that internal homs exist on both sides.

For any object \(x \in\cat{C}\) of a right closed category, a combination of the unit and counit
\begin{equation*}
  \eta^{(x)}_y \colon y \to {[x, x \otimes y]}_r, \qquad \qquad
  \varepsilon^{(x)}_y \colon x \otimes {[x,y]}_r \to y, \qquad \text{for all } y \in\cat{C}
\end{equation*}
gives rise to the morphism
\begin{equation} \label{eq:canonical-morphism-deligne}
  % https://q.uiver.app/#q=WzAsMyxbMCwwLCJbeCwxXV9yXFxvdGltZXMgeSJdLFsyLDAsIlt4LCB4IFxcb3RpbWVzIFt4LDFdX3JcXG90aW1lcyB5XV9yIl0sWzQsMCwiW3gseV0iXSxbMCwxLCJcXGV0YV97W3gsMV1fciBcXG90aW1lcyB5fSJdLFsxLDIsIlt4LCBcXHZhcmVwc2lsb25fMSBcXG90aW1lcyB5XSJdXQ==
  \begin{tikzcd}[ampersand replacement=\&]
    {\phi^{(x)}_y \eqdef {[x,1]}_r\otimes y} \&\& {{[x, x \otimes {[x,1]}_r\otimes y]}_r} \&\& {{[x,y]}_r.}
    \arrow["{\eta^{(x)}_{{[x,1]}_r \otimes y}}", from=1-1, to=1-3]
    \arrow["{[x, \varepsilon^{(x)}_1 \otimes y]}", from=1-3, to=1-5]
  \end{tikzcd}
\end{equation}

The next result, linking the invertibility of these morphisms to a particularly well-behaved type of dualisablility was noted quite early in the development of monoidal categories, see for example~\cite{Kelly1972}.
A concise proof is given in~\cite[Proposition~2.1]{niefield17:coexp}.

\begin{proposition}\label{prop:dualisable-conditions-Deligne}
  The following are equivalent for any object \(x\)
  of a monoidal category \(\cat{C}\):
  \begin{thmlist}[itemsep=1ex]
    \item the internal hom of \(x\) exists, and
    the canonical arrows \(\phi^{(x)}_{y}\) are invertible for all \(y \in\cat{C}\);
    \item the internal hom of \(x\) exists and \(\phi^{(x)}_x\) is an isomorphism; and
    \item there exists an object \(\rd{x} \in\cat{C}\) together with morphisms \(\ev_x\colon x \otimes \rd{x} \to 1\) and \(\coev_x\colon 1 \to \rd{x} \otimes x\), satisfying the snake identities
    \begin{equation}\label{eq:snake-ids}
      \id_x = (\ev_x \otimes \id_x) (\id_x \otimes \coev_x)
      \qquad \text{and} \qquad
      \id_{\rd{x}} = (\id_{\rd{x}} \otimes \ev_x) (\coev_x \otimes \id_{\rd{x}}).
    \end{equation}
  \end{thmlist}
\end{proposition}

We call \(x\) \emph{right (rigidly) dualisable} if any of the equivalent conditions of the previous proposition are met.
In this case, we have \({[x, y]}_r \cong \rd{x} \otimes y\) for all \(y \in \cat{C}\) and in particular \(\rd{x} \cong {[x, 1]}_r\).
The \emph{left (rigid) dualisability} of an object \(x\in\cat{C}\) can be defined similarly by either the invertibility of a canonical morphism
\(\psi^{(x)}_y \colon y \otimes {[x,1]}_{\ell} \to {[x,y]}_{\ell}\) or by the existence of an object \(\ld{x} \in\cat{C}\) endowed with morphisms \(\overline{\ev_x} \colon \ld{x} \otimes x \to 1\) and \(\overline{\coev_x} \colon 1 \to x \otimes \ld{x}\), subject to suitable variants of the snake identities.

\begin{definition}
  A monoidal category \(\cat{C}\) is called \emph{left} or \emph{right rigid} if every object is left or right rigidly dualisable.
  If \(\cat{C}\) is left and right rigid, it is referred to as \emph{rigid}.

  For a rigid category \(\cat{C}\),
  we shall denote its left and right \emph{dualising functors}
  by
  \[
    \ld{(\blank)} \cong {[\blank, 1]}_{\ell}\from \cat{C}^{\op} \to \cat{C},
    \qquad\text{and}\qquad
    \rd{(\blank)} \cong {[\blank, 1]}_{r}\from \cat{C}^{\op} \to \cat{C}.
  \]
\end{definition}

There are categories with only one-sided closedness or rigidity, see~\cite[Theorem~6.3.3]{loregian2021} and~\cite[Example~1.6.2]{Turaev2017}.
Taking the opposite of a rigid category reverses the roles of evaluation and coevaluation. Thus the opposite \(\cat{C}^{\op}\) of the left rigid category \(\cat{C}\) is right rigid.%

The following lemma is well-known, see \eg~\cite[Exercise~2.10.6]{Etingof2015}.

\begin{lemma}\label{lem:transfer-of-rigidity}
  Let \((F, F_2, F_0) \from \cat{C} \to \cat{D}\) be a strong monoidal functor. The image \(Fx\) of  any (rigidly) dualisable object \(x \in\cat{C}\) is dualisable.
\end{lemma}

For example,
if \(x \in \cat{C}\) is an object with right dual \(\rd{x}\),
then \(Fx\) has \(F(\rd{x})\) as its right dual;
the evaluation map is given by
\[
  Fx \otimes F(\rd{x})  \xrightarrow{F_2} F(x \otimes \rd{x}) \xrightarrow{F \ev^{(r)}_x} F1 \xrightarrow{F_0^{-1}} 1,
\]
and the coevaluation is defined analogously.
A straightforward calculation shows the snake identities to be satisfied.

Mapping objects of a (right) rigid monoidal category to their right duals can be extended to a contravariant functor.
This is discussed for example in Section~2.10 of~\cite{Etingof2015} or, in a slightly more general setting, in Section~\ref{sec:tensor-rep->gv}, and leads to the following result.
\begin{theorem}\label{thm:compact-closed->closed}
  Let \(\cat{C}\) be a right rigid category. Then
  \begin{equation}\label{eq:internal hom}
    x \otimes \blank  \adjoint \rd{x} \otimes \blank \qquad \text{for all } x \in \cat{C}.
  \end{equation}
  Furthermore, \(\rd{(\blank)}\) is an equivalence of categories if and only if \(\cat{C}\) is rigid.
\end{theorem}

The main concern of this article is to investigate whether rigidity is indeed equivalent to the tensor-representability of the internal hom, as given in Theorem~\ref{thm:compact-closed->closed}.
To elucidate the underlying problem,
let us assume for a moment that we are given objects \(x, Rx \in \cat{C}\),
such that \(\adj{x \otimes \blank }{Rx \otimes\blank}{\cat{C}}{\cat{C}}\).
One can show that any right rigid dual of \(x\)—if it exists—has to be isomorphic to \(Rx\); \ie, one has \(Rx \cong \rd{x}\) on objects.
Hence, the unit \(\eta^{(x)}_{z} \from z \to Rx \otimes x \otimes z\) and counit \(\varepsilon^{(x)}_{z} \from x \otimes Rx \otimes z \to z \) of the adjunction provide us with natural candidates for the coevaluation and evaluation morphisms:
\begin{equation*}
  \coev_x \eqdef \eta^{(x)}_1 \from 1 \to Rx \otimes x \qquad \qquad \text{and} \qquad \qquad \ev_x \eqdef \varepsilon^{(x)}_1 \from x \otimes Rx \to 1.
\end{equation*}
Evaluating the triangle identities of the adjunction at the monoidal unit yields
\begin{equation*}
  \id_{x} = \varepsilon^{(x)}_{x}  (x \otimes \eta^{(x)}_1) \qquad\qquad\text{and}\qquad\qquad \id_{Rx} = (Rx \otimes \varepsilon^{(x)}_1)  \eta^{(x)}_{Rx}.
\end{equation*}
However, if \(Rx\) is to be a  dual of x in the rigid sense, the snake identities~\eqref{eq:snake-ids} must hold.
For this, we ought to require the stronger condition that
\begin{equation*}
  \varepsilon^{(x)}_x = \varepsilon^{(x)}_1 \otimes \id_x \qquad\qquad\text{and}\qquad\qquad \eta^{(x)}_{Rx} = \eta^{(x)}_1 \otimes Rx.
\end{equation*}
In the next sections we will show that this needs not be true.

\section{Tensor-representability and \GV{} categories}\label{sec:tensor-rep->gv}

In this section we will study closed monoidal categories whose internal homs are tensor-representable;
that is, given by tensoring with dual objects from the left and right, respectively.
Unlike the rigid case, these duals are not required to satisfy the snake identities.
If the  natural transformations relating the left and right sided variants are invertible—see Propositions~\ref{prop:trep-criteria} and~\ref{prop:right-tensor-rep->right-gv}—tensor-representability is a special case of \GV{} duality, which is also called \(*\)-autonomy, see~\cite{barr79,barr95:nonsy,boyarchenko13:groth-verdier}.
In return, we can draw conclusions about the dualising functors.

\begin{definition}\label{def:tensor-representable-endo-functors}
  We call an endofunctor \(F \from\cat{C} \to\cat{C}\) of a monoidal category \emph{right tensor-representable}, if there exists an \(x \in\cat{C}\) such that \(F \cong x \otimes \blank \).
\end{definition}

Theorem~\ref{thm:compact-closed->closed} states that the internal hom of a rigid monoidal category is tensor-re\-pre\-sent\-able in every object.
From now on, for the sake of brevity, we will simply speak of `tensor-representable categories' instead of `closed monoidal categories with tensor-representable internal homs'.
As shown in Lemma~\ref{lem:condition-closed->tensor-dual} below, this is equivalent to the following definition, which does not rely on closedness.

\begin{definition}\label{def:tensor-representable}
  A monoidal category \(\cat{C}\) is called \emph{right tensor-representable}
  if for all \(x \in \cat{C}\) there exists a \emph{right tensor-dual} \(Rx \in\cat{C}\), such that
  \begin{equation} \label{eq:right-ten-dual}
    x \otimes \blank \adjoint Rx \otimes \blank.
  \end{equation}
\end{definition}

Analogous to the previous cases, we call a monoidal category \(\cat{C}\) \emph{left tensor-representable} if \(\cat{C}^{\tensorop}\) is right tensor-representable. A left and right tensor-representable category will be referred to as \emph{tensor-representable}.

Tensor-representable categories encompass all rigid monoidal categories;
as such, some results in the theory of the latter carry over to the former.
For example, the following result is well-known in the case of rigid monoidal abelian categories,
see for example~\cite[Corollary~4.2.13]{Etingof2015}.

\begin{lemma}\label{lem:projectivity-in-tr-cats}
  All objects of an abelian monoidal right tensor-representable category \(\cat{A}\) are projective\footnote{\,%
    In an abelian category \(\cat{A}\),
    an object \(x\) is called \emph{projective}
    if \(\cat{A}(x, \blank)\) is an \emph{exact functor};
    \ie, preserves finite limits and colimits.%
  }
  if and only if its unit is projective.
\end{lemma}
\begin{proof}
  As the first claim implies the second, we only have to show its converse.
  To that end, fix an object \(x \in \cat{A}\).
  The right tensor-representability of \(\cat{A}\) implies that there exists a chain of adjunctions
  \begin{equation*}
    x \otimes \blank \lad Rx \otimes \blank \lad R^2x \otimes \blank \lad \dots,
  \end{equation*}
  which implies \(\cat{A}(x, \blank) \cong \cat{A}(1, Rx \otimes \blank)\).
  The functor \(Rx \otimes \blank \from \cat{A} \to\cat{A}\) has left and right adjoints and is therefore exact.
  Consequently, \(\cat{A}(1, Rx \otimes \blank)\) is—as a composite of exact functors—exact itself.
\end{proof}

While Lemma~\ref{lem:transfer-of-rigidity} holds in the rigid case,
it is not always true for tensor-representable categories.
However, the property is still preserved by monoidal equivalences.

\begin{lemma}\label{lem:transfer-of-tensor-rep}
  Let \(\cat{C}\) be a right tensor-representable category,
  and suppose that \(\cat{D}\) is monoidal.
  Any strong monoidal adjoint equivalence
  \(\adj{F}{F^{-1}}{\cat{C}}{\cat{D}}\)
  equips \(\cat{D}\) with a right tensor-representable structure.
\end{lemma}
\begin{proof}
  Define for any \(x \in \cat{D}\) the object \(\mathsf{R}x \eqdef FRF^{-1}x \in\cat{D}\),
  where \(RF^{-1}x\) is a right tensor-dual of \(F^{-1}x\in\cat{C}\).
  A calculation shows \(\mathsf{R}x \otimes \blank\) to be right adjoint to \(x \otimes \blank\):
  \begin{align*}
    \cat{D}(x \otimes y, z)
    & \cong \cat{C}(F^{-1}(x \otimes y), F^{-1}z)
      \cong \cat{C}(F^{-1}x \otimes F^{-1}y, F^{-1}z)
      \cong \cat{C}(F^{-1}y, RF^{-1}x \otimes F^{-1}z) \\
    & \cong \cat{D}(FF^{-1}y, F(RF^{-1}x \otimes F^{-1}z))
      \cong \cat{D}(y, FRF^{-1}x \otimes z) \eqdef \cat{D}(y, \mathsf{R}x \otimes z).
  \end{align*}
  Thus, \(\cat{D}\) is right tensor-representable.
\end{proof}

Just like in the rigid case,
the relationship between closedness and tensor-repre\-sentability is governed by a family of natural isomorphisms,
as our next result shows.

\begin{lemma}\label{lem:condition-closed->tensor-dual}
  A monoidal category \(\cat{C}\) is right tensor-representable if and only if it is right closed monoidal and for all \(x \in\cat{C}\) there are isomorphisms
  \begin{equation} \label{eq:closed->trep-merely-isos-required}
    \varphi^{(x)}_y \from{[x,y]}_r \to {[x,1]}_r \otimes y \qquad \qquad \text{natural in } y \in\cat{C}.
  \end{equation}
\end{lemma}
\begin{proof}
  If \(\cat{C}\) is right tensor-representable, there exists, by definition, for every \(x \in \cat{C}\) an object \(Rx \in \cat{C}\) such that \(Rx\otimes \blank\) is right adjoint to \(x\otimes \blank\).
  In view of Equation~\eqref{eq:right-internal hom-def}, the category \(\cat{C}\)
  is right closed monoidal with the internal hom at \(x\in \cat{C}\) given by
  \(Rx \otimes \blank\) and the claim follows from the functorial equality  \(Rx \otimes \blank = Rx \otimes 1 \otimes \blank \) for all \(x \in \cat{C}\).

  Conversely, we note that for all objects \(x \in \cat{C}\) we have
  \begin{equation*}
    \cat{C}(x\otimes y, z)
    \cong \cat{C}(y, {[x,z]}_r)
    \cong \cat{C}(y, {[x,1]}_{r}\otimes z) \qquad \text{natural in } y,z\in \cat{C}.\qedhere
  \end{equation*}
\end{proof}

Assume \(\cat{C}\) to be a right tensor-representable category.
By the previous lemma, it is right closed monoidal.
We call \(R \eqdef  {[\blank, 1]}_r \from \cat{C}^{\op} \to \cat{C}\) its \emph{right (tensor-)dualising functor}.
The \emph{left (tensor-)dualising functor} \(L\eqdef {[ \blank , 1]}_\ell\from \cat{C}^{\op}\to \cat{C}\) of a left tensor-representable category is defined analogously.

A right rigid monoidal category is simultaneously left rigid if and only if its right dualising functor is an equivalence of categories.
This is more complicated in the tensor-representable case.
Assume \(\cat{C}\) to be of this type.
Write \(R \from\cat{C}^{\op} \to\cat{C}\) and \(L \from\cat{C}^{\op} \to\cat{C}\) for the right and left tensor-dualising functor, respectively.
Set
\begin{align*}
  \eta^{(x)}_y &\from y \to Rx \otimes x \otimes y, & \varepsilon^{(x)}_y&\from x \otimes Rx \otimes y \to y ,\\
  u^{(x)}_y & \from y \to y \otimes x \otimes Lx, & c^{(x)}_{y} &\from y \otimes Lx \otimes x \to y
\end{align*}
for the units and counits of the corresponding adjunctions.
The left and right tensor dualising functors are related to each other via the natural transformations
\begin{subequations}
  \begin{gather}
    x \xrightarrow{u^{(Rx)}_x} x \otimes Rx \otimes LRx \xrightarrow{\varepsilon^{(x)}_1 \otimes LRx}LRx, \label{eq:can-morph-unit} \\
    x \xrightarrow{\eta^{(Lx)}_x} RLx \otimes Lx \otimes x \xrightarrow{RLx \otimes c^{(x)}_1}RLx. \label{eq:can-morph-counit}
  \end{gather}
\end{subequations}

\begin{proposition}\label{prop:trep-criteria}
  Let \(\cat{C}\) be tensor-representable.
  The right dualising functor \(R\from\cat{C}^{\op} \to\cat{C}\) is right adjoint to \(L^{\op}\from\cat{C}\to\cat{C}^{\op}\).
  It is an equivalence  if and only if the canonical morphisms of Equations~\eqref{eq:can-morph-unit} and~\eqref{eq:can-morph-counit} are invertible.
  In this case \(L^{\op}\) is a quasi-inverse of \(R\).
\end{proposition}
\begin{proof}
  For any two objects \(x,y \in\cat{C}\), we have
  \begin{equation*}
    \cat{C}^{\op}(Lx, y) =  \cat{C}(y, Lx) \cong  \cat{C} (y \otimes x, 1) \cong  \cat{C}(x,  Ry).
  \end{equation*}
  A direct computation shows that the unit and counit of this adjunction are given by the natural transformations displayed in Equation~\eqref{eq:can-morph-unit} and~\eqref{eq:can-morph-counit}.
  If these morphisms are invertible, \(R\) and \(L^{\op}\) are quasi-inverse to each other.

  Conversely, assume \(R\) to be an equivalence of categories.
  Since any equivalence can be bettered to an adjoint equivalence and adjoints are unique up to unique natural isomorphism, \(L^{\op}\) must be a quasi-inverse to \(R\) and the unit and counit morphisms stated in Equation~\eqref{eq:can-morph-unit} and~\eqref{eq:can-morph-counit} are invertible.
\end{proof}

Monoidal categories endowed with a type of duality which is governed by an anti-auto-equivalence were introduced and studied under the name of \(*\)-autonomous categories by Barr, see~\cite{barr79,barr95:nonsy}.
Boyarchenko and Drinfeld, who referred to them as \GV{} categories, see~\cite{boyarchenko13:groth-verdier}, investigated their applications and properties in various algebraic contexts.
We are going to loosely follow their notation.

\begin{definition}\label{def:*-autonomous-drinfeld}
  A \emph{(right) \GV{} category} comprises a pair \((\cat{C}, d)\) of a monoidal category \(\cat{C}\) and an
  object \(d\in \cat{C}\), such that for all \(x\in\cat{C}\) there exists an object \(D_r x \in \cat{C}\) satisfying
  \begin{equation} \label{eq:*-aut-condition}
    \cat{C}(x \otimes \blank ,d) \cong\cat{C}(\blank, D_{r} x),
  \end{equation}
  and the induced \emph{right dualising functor} \(D_r\from \cat{C}^{\op}\to\cat{C}\) is an equivalence of categories.

  If \(d=1\) is the monoidal unit, one speaks of an \emph{\(r\)-category}.
\end{definition}

Because \(D_r\) is an equivalence,
one also obtains a natural isomorphism for its inverse:
\begin{equation} \label{eq:gv-inverse-condition}
  \cat{C}(x \otimes y, d) \cong \cat{C}(y, D_r x) \cong \cat{C}(x, D_r^{-1} y), \qquad \text{for all } x, y \in \cat{C}.
\end{equation}

\begin{remark}\label{rmk:right-gv-vs-the-world}
  In~\cite{boyarchenko13:groth-verdier}, the authors instead consider
  a pair \((\cat{C}, d)\)
  subject to the natural isomorphism \(\cat{C}(\blank \otimes y, d) \cong\cat{C}(\blank, D_{\ell} y)\),
  such that the induced functor \(D_{\ell}\from \cat{C}^{\op} \to \cat{C}\) is an equivalence.
  In our framework, we shall call this a \emph{left} \GV{} category structure on \(\cat{C}\),
  with \emph{left} dualising functor \(D_{\ell}\).
  Equations~\eqref{eq:*-aut-condition}
  and~\eqref{eq:gv-inverse-condition} become
  \begin{equation}\label{eq:left-gv-condition-and-inverse}
    \cat{C}(x \otimes y, d) \cong \cat{C}(x, D_{\ell}y) \cong \cat{C}(y, D^{-1}_{\ell}x).
  \end{equation}
\end{remark}

\begin{remark}\label{rmk:GV-is-symmetric}
  Because the dualising functor of a \GV{} category is assumed to be an equivalence,
  \GV{} structures are always two-sided:
  in view of Equation~\ref{eq:gv-inverse-condition},
  a right \GV{} category \((\cat{C}, d)\)
  is also a left \GV{} category
  with the same dualising object \(d\), as well as \(D_r^{-1} \eqdef D_{\ell}\).
  Likewise, all left \GV{} categories are also right-sided.

  However,
  for the sake of clarity and to denote the directional `bias',
  \GV{} structures will nevertheless often explicitly carry a prefix.
  In particular, in Section~\ref{sec:tensor-rep-functor-cats}
  we consider the \GV{} structure of the category of copresheaves \(\Cpsh\) over a \GV{} `base' category \(\cat{C}\).
  In this case, the \emph{left} dualising functor of \(\cat{C}\) naturally gives rise to a \emph{right} dualising functor for \(\Cpsh\); see Proposition~\ref{prop:fun-cat-GV}.
\end{remark}

The following Proposition, also discussed in Section~2 of~\cite{boyarchenko13:groth-verdier}, is a direct consequence of Equations~\eqref{eq:*-aut-condition} and~\eqref{eq:gv-inverse-condition}.

\begin{proposition}\label{prop:gv->closed}
  Every (right) \GV{} category \((\cat{C}, d)\) is closed monoidal;
  for all \(x, y \in \cat{C}\), the left and right internal homs are given by
  \begin{equation} \label{eq:gv->internal hom}
    \begin{aligned}
      {[x, y]}_{\ell} &\eqdef D_r^{-1}(x \otimes D_r y) \eqdef D_{\ell}(x \otimes D_{\ell}^{-1} y), \\
      {[x, y]}_r &\eqdef D_r(D_r^{-1}y \otimes x) \eqdef D_{\ell}^{-1}(D_{\ell}y \otimes x).
    \end{aligned}
  \end{equation}
\end{proposition}

One can use Equation~\eqref{eq:gv->internal hom} to recover the dualising functor of a \GV{} category
from the relevant internal hom by evaluating it at the dualising object in its second component.

\begin{remark}\label{rmk:gv-d-from-D}
  If \((\cat{C},d)\) is a (left or right) \GV{} category,
  one can recover its dualising object by
  applying the dualising functor to the monoidal unit of \(\cat{C}\).
  In fact, as shown in~\cite[Remark~2.1.(4)]{boyarchenko13:groth-verdier}, this relationship is bidirectional:
  \begin{equation}\label{eq:gv:1:from:d:and:d:from:1}
    D_r 1 \cong d, \qquad\!
    D_r d \cong 1, \qquad\!
    D_{\ell} 1 \cong d, \qquad\!
    D_{\ell} d \cong 1, \qquad\!
    D_r^2 1 \cong 1, \qquad\!
    D_r^2d \cong d.
  \end{equation}

  In combination with Proposition~\ref{prop:gv->closed},
  one obtains an explicit description of the dualising functors of \(\cat{C}\):
  \begin{equation}\label{eq:gv-dualising-from-internal-hom-vice-versa}
    D_r x \cong {[x, d]}_r \qquad \text{ and } \qquad
    D_{\ell}x \eqdef D_r^{-1}x \cong {[y, d]}_{\ell}, \qquad \text{for all }x \in \cat{C}.
  \end{equation}
\end{remark}

The `symmetric' nature of \GV{} categories is in stark contrast to the one-sided notions of duality such as rigidity and tensor-reprensentability.

\begin{example}\label{ex:right-rigid-not-gv}
  Consider a Hopf algebra \(H\) whose antipode is not invertible as constructed for example in \cite{schauenburg2000:FaithfulFlatnessHopf}.
  Its finite-dimensional right comodules are a left-rigid category \(\rcomod{H}\) which as discussed in Remark~2 of \cite{ulbrich1990:HopfAlgebrasRigid} cannot be right rigid.
  If \(\rcomod{H}\) were a (left) \GV{} category, there would exist, as a consequence of Remark~\ref{rmk:gv-d-from-D}, a comodule \(M \in \rcomod{H}\) such that \(\ld{( \blank )}\otimes M \from \rcomod{H}^{\op}\to \rcomod{H}\) is an equivalence of categories.
  This implies that \(M\) is one-dimensional and therefore the left dualising functor \( \ld{( \blank )}\from \rcomod{H}^{\op}\to \rcomod{H}\) would be an equivalence which is a contradiction.
\end{example}

An object \(\alpha\in \cat{C}\) is called \emph{invertible} if there exists a \(\beta\in \cat{C}\) such that \(\alpha \otimes \beta \cong 1 \cong \beta \otimes \alpha\).
The inverse \(\beta\) is unique up to isomorphism and we write \(\alpha^{-1}\eqdef \beta\).
Proposition~2.3 of \cite{boyarchenko13:groth-verdier} shows that invertible objects govern the possible choices for dualising objects of a \GV{} category.

\begin{lemma}\label{lem:choices-of-dualising-objects}
  Given a right \GV{} category \((\cat{C},d)\), there is an anti-equiva\-lence between its full subcategories of invertible and dualising objects given by
  \begin{equation*}
    \mathrm{Inv}(\cat{C}) \to \mathrm{Dual}(\cat{C}), \qquad\alpha \mapsto \alpha^{-1}\otimes d.
  \end{equation*}
\end{lemma}

One may also use~\cite[Proposition~2.3]{boyarchenko13:groth-verdier} to produce examples of closed categories that cannot be equipped with any \GV{} dualising object.

\begin{example}\label{rmk:closed-but-not-gv}
  The category \((\mathsf{Set}, \times, \1)\) of sets
  with its usual Cartesian monoidal structure is closed monoidal.
  If some set \(X\) is invertible in the above discussed sense
  then it must already be isomorphic to the terminal set \(\1\).
  Thus, it is enough to see that \(\mathsf{Set}(\blank, \1)\from \mathsf{Set}^{\op} \to \mathsf{Set}\) is not an anti-equivalence,
  which follows by \(\1\) being terminal, to deduce that there is no \GV{} structure on this category.
\end{example}

Let us now characterise the relationship between \GV{} and tensor-representable categories.

\begin{proposition}\label{prop:right-tensor-rep->right-gv}
  For a tensor-representable category  \(\cat{C}\) the following are equivalent:
  \begin{thmlist*}
    \item the right dualising functor admits a quasi-inverse, \label{itm: dualising-functor-trep-invertible}
    \item \(\cat{C}\) is an \(r\)-category, and \label{itm: trep-is-r-cat}
    \item there is an object \(d\in \cat{C}\) such that \((\cat{C}, d)\) is a right \GV category. \label{itm: trep-is-GV}
  \end{thmlist*}
\end{proposition}
\begin{proof}
  By definition we have \(\ref{itm: dualising-functor-trep-invertible} \implies \ref{itm: trep-is-r-cat} \implies \ref{itm: trep-is-GV}\).

  Therefore, let us assume that \(\cat{C}\) is \GV{} with respect to the dualising object \(d\in \cat{C}\).
  We write \(L, R \from \cat{C}^{\op} \to \cat{C}\) for its left and right tensor-dualising functors and
  \(D_r\from \cat{C}^{\op} \to \cat{C}\) for the dualising functor induced by \(d\in \cat{C}\).
  We have
  \begin{align*}
    Rd \otimes d
    \overset{\eqref{eq:closed->trep-merely-isos-required}}{\cong} [d,d]_r
    & \overset{~\eqref{eq:gv->internal hom}}{\cong} D_r(D_r^{-1} d \otimes d)
      \overset{~\eqref{eq:gv:1:from:d:and:d:from:1}}{\cong} D_rd
      \overset{~\eqref{eq:gv:1:from:d:and:d:from:1}}{\cong} 1,\\
    d \otimes Ld
    \overset{\eqref{eq:closed->trep-merely-isos-required}}{\cong} [d,d]_{\ell}
    & \overset{~\eqref{eq:gv->internal hom}}{\cong} D_r^{-1}(d \otimes D_r d)
      \overset{~\eqref{eq:gv:1:from:d:and:d:from:1}}{\cong} D_r^{-1} d
      \overset{~\eqref{eq:gv:1:from:d:and:d:from:1}}{\cong} 1
  \end{align*}
  and therefore \(d\otimes Rd \cong d \otimes Rd \otimes d \otimes Ld \cong 1\).
  In other words, \(d\) is invertible and \(1 \cong Rd\otimes d\) another dualising object due to Lemma~\ref{lem:choices-of-dualising-objects}.
  Thus, \(\ref{itm: trep-is-GV} \implies \ref{itm: trep-is-r-cat}\).

  To complete the proof it suffices to observe that the right dualising functor of \(\cat{C}\) considered as an \(r\)-category and its right tensor-dualising functor coincide.
\end{proof}

Conversely, it is not true that an \(r\)-category is tensor-representable.

\begin{example}\label{ex: r-cat-not-always-trep}
  Let \(A\) be a commutative Frobenius algebra over a field \(\k\).
  We write \(\lmod{A}\) for the category of \emph{finitely-generated} left \(A\)-modules.
  One can show, see for example~\cite[Sections 15 and 16]{Lam1999}, that for all \(M \in \lmod{A}\) the canonical morphism
  \begin{equation*}
    \phi_M \from M \to  \rrd M\eqdef \Hom_A(\Hom_A(M,A),A) , \qquad \qquad \phi_M(x)\alpha = \alpha(x)
  \end{equation*}
  is an isomorphism.
  Put differently, \(\lmod{A}\) is an \(r\)-category with \(\Hom_A(\blank,A) \from \cat{C}^{\op} \to\cat{C}\) as dualising functor.
  A classical result in representation theory shows that
  for any commutative algebra \(A\) and \(A\)-module \(M\) there exists an \(A\)-module \(N\)
  such that \(M \otimes \blank \lad N \otimes \blank\) if and only if \(M\) is finitely-generated projective.
  For a concise proof, we refer the reader to~\cite[Proposition~2.1]{niefield17:coexp} or Section~\ref{sec:tensor-rep-functor-cats}.
  Since the finitely-generated projective \(A\)-modules are (rigidly) dualisable, we obtain that \(\lmod{A}\) is tensor-representable if and only if it is rigid if and only if all (finitely-generated) \(A\)-modules are projective, which in turn is equivalent to \(A\) being semisimple.

  However, not all commutative algebras and in particular Frobenius algebras are semisimple.
  Consider for example the truncated polynomial ring \(\k [x]/(x^{2})\).
  It is not semisimple since its regular module is indecomposable but contains the one-dimensional simple submodule generated by \(x\).
  Moreover, the functional
  \begin{equation*}
    \phi \from \k[x]/(x^{2}) \to \k, \qquad ax + b \mapsto a
  \end{equation*}
  turns it into a Frobenius algebra.
\end{example}

To conclude this section, we will abstractly compare tensor-representable \(r\)-categories and rigid monoidal categories.

Let \(\cat{C}\) be a right tensor-representable category whose dualising functor \(R \from\cat{C}^{\op} \to\cat{C}\) is an equivalence.
For all \(x,y \in\cat{C}\), there is a canonical isomorphism
\[
  \tau_{x,y} \from\cat{C}(y\otimes x \otimes Rx \otimes Ry,1) \to\cat{C}^{\op}(y\otimes x, R^{-1}(Rx \otimes Ry)),
\]
see~\cite[Section~4]{boyarchenko13:groth-verdier}.
It is given by the identifications
\begin{equation*}
  \cat{C}(y\otimes x\otimes Rx \otimes Ry, 1) \cong \cat{C}(Rx \otimes Ry, R(y \otimes x)) \cong\cat{C}^{\op}(y \otimes x, R^{-1}(Rx \otimes Ry)).
\end{equation*}
As shown in Lemma 4.1 and Remark 4.2 of~\cite{boyarchenko13:groth-verdier}, the natural transformation
\begin{equation*}
  \mu_{x,y} = R(\tau_{x,y}(\varepsilon^{(y)}_1(y \otimes \varepsilon^{(x)}_1 \otimes Ry)))\from Rx \otimes Ry \to R(x \otimes^{\op} y) \eqdef R(y \otimes x)
\end{equation*}
endows \(R\from\cat{C}^{\op,\tensorop} \to\cat{C}\) with the structure of a lax monoidal functor.
That is, for all \(x,y,z \in\cat{C}\) we have
\begin{equation*}
  \mu_{x, y\tensorop z}(Rx \otimes \mu_{y,z}) = \mu_{x \tensorop y, z}(\mu_{x,y}\otimes Rz),
  \qquad
  \mu_{1,x}(R1 \otimes Rx) = \id_{Rx}= \mu_{x,1}(Rx \otimes R1).
\end{equation*}
As a consequence of~\cite[Proposition~4.4]{boyarchenko13:groth-verdier}, we get the following characterisation of rigidity.
\begin{proposition}\label{prop:GV-rigid-characterisation}
  Let \(\cat{C}\) be a right tensor-representable category whose dualising functor \(R \from\cat{C}^{\op, \tensorop} \to\cat{C}\) is an equivalence.
  It is rigid if and only if \((R, \mu,R1)\) is strong monoidal; \ie,
  \(\mu_{x,y} \from Rx \otimes Ry \to R(x \tensorop y)\) is an isomorphism for all \(x, y\in\cat{C}\).
\end{proposition}

\section{Tensor-representability for functor categories}\label{sec:tensor-rep-functor-cats}

We will now investigate the previously discussed types of dualities in the context of representation theory.
The starting point is the convolution closed monoidal structure on suitable subcategories of functor categories;
examples of which include \(\mathfrak{sl}_{2}(\mathbb{C})\)-crystals, modules over certain path algebras, Mackey functors, and
group-graded representations arising from crossed modules.
Before examining them in Section~\ref{sec:mky-2-grps}, we develop abstract tools to detect \GV{} structures, tensor-representability, and rigidity—see Proposition~\ref{prop:fun-cat-GV} and its Corollary~\ref{cor:day-compact}.
We conclude the section with Proposition~\ref{prop:characterisation-of-rigidity}, in which tools are developed to characterise when rigidity and tensor-representability for functor categories coincide.
\bigskip

Throughout this section, fix a field \(\k\) and let \(\kvect\) as well as \(\kVect\) denote strict monoidal categories equivalent to the categories of (finite-dimensional) \(\k\)-vector spaces.
We write \(\rd*{(\blank)}\from \kvect^{\op} \to \kvect\) for the \(\k\)-linear dualising functor.
Unless explicitly stated, from now on all (monoidal) categories and functors will be enriched over \(\kVect\).

\begin{hypothesis}\label{hypo:eso-enriched}
  In the following \(\cat{C}\) will be an essentially small \(\k\)-linear category.
\end{hypothesis}
For the sake of brevity, the prefixes `enriched', `\(\k\)-linear', and `(essentially) small' will often be omitted.
The category of functors from \(\cat{C}\) to \(\kVect\) will be denoted by \(\Cpsh\).

\subsection{Linear functor categories}\label{sec:lin-fun-cats}

In order to define an appropriate monoidal structure for functors,
one needs the notion of an enriched coend.
We will briefly recall it in this section,
specialising the definition for \(\kVect\).
More general accounts can be found for example
in~\cite[Section~2.3]{gray80:closed},~\cite[Chapter~3]{kelly05:basic},
and~\cite[Chapter~4]{loregian2021}.

\begin{definition}\label{def:coend-vect}
  Let \(\cat{C}\) be an essentially small \(\k\)-linear category,
  and consider a \(\k\)-linear functor \(P\from \cat{C}^{\op} \otimes_{\k} \cat{C} \to \kVect\).
  The \emph{coend} of \(P\) is the coequaliser in \(\kVect\):
  \begin{equation}
    \begin{tikzcd}
      \displaystyle{\coprod_{a, b \in \cat{C}} \cat{C}(a, b) \otimes_{\k} P(b, a)} & \displaystyle{\coprod_{a \in \cat{C}} P(a,a)} & \displaystyle{\int^{a \in \cat{C}} P(a,a).}
      \arrow[shift right=1, from=1-1, to=1-2]
      \arrow[shift left=1, from=1-1, to=1-2]
      \arrow[two heads, from=1-2, to=1-3]
    \end{tikzcd}
  \end{equation}
  For \(a, b \in \cat{C}\) the two parallel morphisms are given by
  \begin{align*}
    \cat{C}(a,b)\otimes_{\k} P(b, a) &\to P(a,a), \qquad\qquad  f \otimes x \mapsto P(f, a)x, \\
    \cat{C}(a, b) \otimes_{\k} P(b, a) &\to P(b, b), \qquad\qquad g \otimes y \mapsto P(b,g)y.
  \end{align*}

  The \emph{end} \(\int_{a \in\cat{C}} P(a,a)\) of \(P\) may be defined analogously as an equaliser in \(\kVect\).
\end{definition}
A consequence of \(\cat{C}\) being essentially small is that we can rewrite the (co)equalisers defining the respective (co)ends in Definition~\ref{def:coend-vect} to index over a set of objects, instead of a proper class.
In particular, since \(\kVect\) is complete and cocomplete, the end and coend of \(P\from \cat{C}^{\op} \otimes_{\k} \cat{C} \to \kVect\) exists.

In case the indexing category \(\cat{C}\) can be inferred unambiguously from the context, we will omit writing it explicitly.

Ends and coends have a number of useful properties, two of which we are going to frequently use in the following.
First, they are functorial: Any natural transformation \(\mu \from P \to P'\) between functors \(P, P' \from\cat{C}^{\op} \times\cat{C} \to \kVect\) induces via the universal property of (co-)equalisers a unique morphism \(\int\! P \xrightarrow{\;\int\! \mu\;} \int\! P'\).

Second, for appropriate functors \(P\),
the \emph{Fubini–Tonelli interchange law} holds:
if any of
\begin{equation*}
  \int^{a \in \cat{C}}\!\!\!\! \int^{b \in \cat{C}}\!\! P(a, a, b, b), \qquad
  \int^{(a, b) \in \cat{C} \times \cat{C}}\!\! P(a, a, b, b), \qquad
  \int^{b \in \cat{C}}\!\!\!\! \int^{a \in \cat{C}}\!\! P(a, a, b, b)
\end{equation*}
exist, they all do and are isomorphic.
We will denote any of the three isomorphic terms above as \(\int^{a,b} P(a,a,b,b)\).

\begin{lemma}[enriched Yoneda]\label{lem:enriched-yoneda}
  Let \(\cat{C}\) be a category,
  and suppose that \(F \from \cat{C} \to \kVect\) and \(G \from \cat{C}^{\op} \to \kVect\)
  are functors.
  There are natural isomorphisms
  \begin{equation}\label{eq:coYoneda-covariant}
    F \cong \int^{c} \cat{C}(c, \blank) \otimes_{\k} Fc \cong \int_{c} \kVect(\cat{C}(\blank, c), Fc),
  \end{equation}
  and
  \begin{equation}\label{eq:coYoneda-contravariant}
    G \cong \int^{c} \cat{C}(\blank, c) \otimes_{\k} Gc \cong \int_{c} \kVect(\cat{C}(c, \blank), Gc).
  \end{equation}
\end{lemma}
As the name suggests,
the above lemma is an analogue of the Yoneda lemma in the enriched case;
as such, we will only refer to it as `the Yoneda lemma' in subsequent considerations.
For a proof, see~\cite[Section~2.4]{kelly05:basic}.

\begin{definition}\label{def:day-convolution}
  Let \(\cat{C}\) be a monoidal category.
  Given two \(\k\)-linear functors \(F, G \in \Cpsh\),
  the \emph{Day convolution} of \(F\) and \(G\) is given by the coend
  \begin{equation} \label{eq:day-convolution}
    F \ostar G \eqdef \int^{a, b} \cat{C}(a \otimes b, \blank) \otimes_{\k} Fa \otimes_{\k} Gb.
  \end{equation}
\end{definition}

Assuming the monoidal category \(\cat{C}\) is left or right closed monoidal,
we may simplify the Day convolution product using the Yoneda lemma for coends:
\begin{subequations}
  \begin{align}
    \begin{aligned}
      F \ostar G &= \int^{a,b }\cat{C}(a \otimes b, \blank) \otimes_{\k} Fa \otimes_{\k} Gb
                   \cong \int^{a,b }\cat{C}(a, {[b, \blank]}_{\ell}) \otimes_{\k} Fa \otimes_{\k} Gb \\
                 & \overset{\eqref{eq:coYoneda-covariant}}{\cong} \int^{b } F({[b,\blank]}_{\ell}) \otimes_{\k} Gb,
    \end{aligned}\label{eq:day-convolution-simp-left}\\
    \begin{aligned}
      F \ostar G &= \int^{a,b }\cat{C}(a \otimes b, \blank) \otimes_{\k} Fa \otimes_{\k} Gb
                   \cong \int^{a,b }\cat{C}(b, {[a, \blank]}_r) \otimes_{\k} Fa \otimes_{\k} Gb \\
                 & \overset{\eqref{eq:coYoneda-covariant}}{\cong} \int^{a } Fa \otimes_{\k} G({[a,\blank]}_r).
    \end{aligned}\label{eq:day-convolution-simp-right}
  \end{align}
\end{subequations}

An illuminating example explaining the convolution tensor product is given by considering a \(\k\)-linear monoidal category \(\cat{X}\) with a single object \(x\).
We may identify \(\cat{X}\) with the commutative algebra \(A \eqdef \End_{\cat{X}}(x)\) and  \([\cat{X}, \kVect] \cong \rMod{A}\).
Now consider two functors \(F, G \from\cat{X} \to \kVect\) and write \(M \eqdef Fx,\ N\eqdef Gx\) for their corresponding modules over \(A\).
By Equation~\eqref{eq:day-convolution-simp-left} and the definition of coends,
\begin{equation*}
  (F \ostar G)x \cong M \otimes_{\k}  N / \langle ma \otimes_{\k} n - m\otimes_{\k} an \mid m \in M, n \in N, a \in A\rangle = M \otimes_{\!A} N,
\end{equation*}
where \(an\eqdef na\) for all \(a\in A\) and \(n\in N\).
In other words, one recovers the tensor product of modules over commutative algebras.

\begin{theorem}[{\cite{day71:const}}]\label{thm:closed-mon-str-on-functors}
  For any monoidal category \(\cat{C}\),
  the category \(\Cpsh\) is closed monoidal
  with Day convolution as its tensor product,
  \(\cat{C}(1, \blank)\) as its unit,
  and the internal homs given for all \(F, G \in \Cpsh\) by
  \begin{subequations}
    \begin{align}
      {[F, G]}_{\ell} &\eqdef \int_{a, b} \kVect(\cat{C}(\blank \otimes a, b), \kVect(Fa, Gb)), \label{eq:Day-internal-hom-left}\\
      {[F, G]}_r &\eqdef \int_{a, b} \kVect(\cat{C}(a \otimes \blank, b), \kVect(Fa, Gb)). \label{eq:Day-internal-hom-right}
    \end{align}
  \end{subequations}
\end{theorem}

For the rest of this article,
the category of functors \(\Cpsh\) of a monoidal category \(\cat{C}\) to \(\kVect\)
shall always be equipped with the closed monoidal structure given by Theorem~\ref{thm:closed-mon-str-on-functors}.

If the monoidal category \(\cat{C}\) is itself closed,
then the formulas for the internal homs of \(\Cpsh\) may be simplified
by means of the Yoneda lemma:
\begin{gather}
  \begin{aligned} \label{eq:right-internal hom-simplification}
    &{[F,G]}_r
      = \int_{a,b} \kVect (\cat{C}(a \otimes \blank, b), \kVect(Fa, Gb)) \\
    &\quad \cong \int_{a, b} \kVect(\cat{C}(a \otimes \blank , b) \otimes_{\k}Fa, Gb)
      \cong \int_b \kVect\Big(\!\int^a\!\!\cat{C}(a \otimes \blank , b) \otimes_{\k}Fa , Gb\Big)\\
    &\quad \cong \int_b \kVect\Big(\!\int^a\!\!\!\cat{C}(a   ,{[\blank, b]}_{\ell}) \otimes_{\k}Fa , Gb\Big)
      \overset{\eqref{eq:coYoneda-covariant}}{\cong} \int_b \kVect(F{[\blank,b]}_{\ell}, Gb),
  \end{aligned} \\
  \begin{aligned} \label{eq:left-internal hom-simplification}
    {[F,G]}_{\ell}
    \cong \int_b \kVect\Big(\int_a \cat{C}(a, {[\blank, b]}_r) \otimes_{\k}Fa , Gb\Big)
    \overset{\eqref{eq:coYoneda-covariant}}{\cong} \int_a \kVect(F{[\blank,b ]}_r, Gb).
  \end{aligned}
\end{gather}

For later applications, we remark that \(\kVect\) naturally acts on \(\Cpsh\), see \cite[Section~3.7]{kelly05:basic}.
\begin{definition}\label{def:copowering}
  A \(\kVect\)-enriched category \(\cat{C}\) is \emph{copowered} over \(\kVect\) if there exists for each pair \(c\in \cat{C}\) and \(v \in \kVect\) an object \(c \boxtimes v\) such that
  \begin{equation}\label{eq: copowering}
    \cat{C}(c\boxtimes v, d) \cong \kVect(v, \cat{C}(c,d)), \qquad\quad \text{natural in } c,d \in\cat{C}, v \in \kVect.
  \end{equation}
\end{definition}
A variation of the parameter theorem, see~\cite[Section IV.7]{MacLane1998}, shows that the object assignment \(- \boxtimes \bblank \from \cat{C} \otimes_{\k} \kVect \to \cat{C}\) lifts to a functor.
Moreover for every \(c\in \cat{C}\), Equation~\eqref{eq: copowering} states that \(c \boxtimes -\) is the left adjoint to the Hom functor \( \cat{C}( c , -)\).
Therefore any copowering of \( \cat{C}\) is unique up to unique isomorphism.
\begin{lemma}\label{lem:copowering-presheaves}
  The presheaves \(\Cpsh\) of a \(\kVect\)-enriched category \( \cat{C}\) are copowered over \(\kVect\) via
  \begin{equation*}
    \blank \boxtimes \bblank \from \Cpsh \otimes_{\k} \kVect \to \Cpsh, \qquad (F\boxtimes v)c = Fc \otimes_{\k} v
  \end{equation*}
  and for all \(F,G \in \Cpsh\) and \(v \in \kVect\) we have
  \begin{equation}\label{eq:compat-Day-convolution-copowering}
    (F \ostar G) \boxtimes v = F \ostar (G \boxtimes v)
  \end{equation}
\end{lemma}
\begin{proof}
  We fix \(F, G \in \Cpsh\) and \(v\in V\).
  Using the description of the morphism spaces of \(\Cpsh\) given in \cite[Chapter~2]{kelly05:basic}, we compute
  \begin{align*}
    \Cpsh(F\boxtimes v, G)
    & \cong \int_{c} \kVect(Fc \otimes_{\k} v, Gc)
      \cong \int_{c} \kVect(v, \kVect(Fc, Gc)) \\
    & \cong \kVect(v, \int_{c} \kVect(Fc, Gc))
      \cong \kVect(v, \Cpsh(F,G)).
  \end{align*}
  To prove the second claim, we simply observe that
  \begin{equation*}
    (F\ostar G)\boxtimes v = \int^{a,b} \cat{C}(a \otimes b, \blank) \otimes_ {\k} Fa \otimes_{\k} Gb\otimes_{\k} v
    =F \ostar(G\boxtimes v).\qedhere
  \end{equation*}
\end{proof}

For any vector space \(V \neq \{0\}\) and \(W\in \kVect\) the canonical map
\begin{equation*}
  \eta^{(V)}_{W}\from W \to \kVect(V,V\otimes_{\k} W), \qquad w \mapsto \blank \otimes_{\Bbbk} w
\end{equation*}
is injective.
Since \(\eta^{(V)} \from \Id \to \kVect(V, V\otimes_{\Bbbk} \blank )\) is the unit of the tensor-hom adjunction \(\adj{V\otimes_{\k} \blank }{\kVect(V, \blank)}{\kVect}{\kVect}\), it follows that \(V\otimes_{\k} \blank \) is faithful.
This readily lifts to functor categories.
\begin{corollary}\label{cor:copwering-is-faithful}
  For each \(0 \neq F\in \Cpsh\) the functor \(F \boxtimes \blank \from \kVect \to \Cpsh\) is faithful.
\end{corollary}
\begin{proof}
  Since \(F\neq 0\) there exists a \(c\in \cat{C}\) such that \(Fc \neq 0\) and the claim follows from the faithfulness of \(Fc \otimes_{\k} \blank \).
\end{proof}

\subsection{Detecting \GV{} structures and tensor-representability}\label{sec:tens-repr-mack}

We will now see that if the duality structures of the source and target categories are compatible in a certain way,
the corresponding functor categories are either \GV{} or—if stronger conditions are met—even tensor-representable.
Inspired by arguments of Day, see~\cite{day06:compac}, we prove that in the semisimple case these conditions are automatically satisfied.

\begin{definition}\label{def:hom-finite-category}
  A category \(\cat{C}\) is called \emph{hom-finite},
  if \(\dim(\cat{C}(a,b))< \infty\) for all \(a, b \in\cat{C}\).
\end{definition}

For a hom-finite category \(\cat{C}\),
instead of studying functors \(F\from \cat{C} \to \kVect\) to \emph{all} vector spaces,
we are interested in the full subcategory \(\Cpsh_{\fin}\) of (point-wise)
\emph{finite-dimensional} functors of \(\Cpsh\).
That is, \(F \in \Cpsh\) such that \(Fx\) is finite-dimensional for all \(x \in \cat{C}\).
\begin{hypothesis}\label{hypothesis:tens-repr-mack}
  Within Section~\ref{sec:tens-repr-mack}, we assume that all relevant ends and coends are finite-dimensional vector spaces.
\end{hypothesis}

\begin{remark}\label{rmk:finiteness-implies-hypothesis}
  If \( \cat{C}\) has finitely many objects and only finite-dimensional morphism spaces, the end and coend of any suitable functor \(F \from \cat{D}_{1}\otimes_{\k} \dots \otimes_{\k}\cat{D}_{n} \to \kvect\), where \( \cat{D}_{i} \in\{ \cat{C}, \cat{C}^{\op}\}\) exists and is finite-dimensional.
  In particular, Hypothesis~\ref{hypothesis:tens-repr-mack} holds.
\end{remark}

Hypothesis~\ref{hypothesis:tens-repr-mack}
implies that \(\Cpsh_{\fin}\) is closed monoidal;
that is, the tensor product and the internal hom are finite-dimensional at every point.
These assumptions impose a certain finiteness condition on \(\cat{C}\) itself, or—as the next example will show—on a dense\footnote{\,%
  A subcategory \(\cat{O}\) of \(\cat{C}\) is \emph{dense} if
  the restricted Yoneda embedding \(\cat{C} \to [\cat{O}^{\op}, \kVect]\)
  is fully faithful.%
}
subcategory of it.

\begin{example}\label{ex:mackey-dense-subcategories}
  Let \(\cat{C}\) be a hom-finite category,
  and assume that there is a full subcategory \(\cat{S}\) such that every object of \(\cat{C}\) may be written as a finite direct sum of objects in \(\cat{S}\).
  In this setting, one can reformulate Lemma~\ref{lem:enriched-yoneda} to index over this finite set:
  for example, given \(F\from \cat{C}\to \kvect\), we have
  \begin{equation} \label{eq:reduction-via-dense}
    \begin{aligned}
      Fx &\cong F\Big(\!\bigoplus_{i = 1}^n s_i\Big)
           \cong \bigoplus_{i = 1}^n Fs_i
           \cong \bigoplus_{i = 1}^n \int^{s \in \cat{S}}\!\! \cat{S}(s, s_i) \otimes Fs
           \cong \bigoplus_{i = 1}^n \int^{s \in \cat{S}}\!\! \cat{C}(s, s_i) \otimes Fs \\
         &\cong \int^{s \in \cat{S}}\!\! \cat{C}\Big(s, \bigoplus_{i = 1}^n s_i\Big) \otimes Fs
           \cong \int^{s \in \cat{S}}\!\! \cat{C}(s, x) \otimes Fs.
    \end{aligned}
  \end{equation}
  Panchadcharam and Street used this `reduction' to show that finite-dimensional Mackey functors are a \GV{} category, see~\cite[Section~9]{panchadcharam07:mackey}.
\end{example}

Hypothesis~\ref{hypothesis:tens-repr-mack}
does not only imply a closed structure on \(\Cpsh_{\fin}\),
but yields a canonical notion of a dual.

\begin{proposition}\label{prop:fun-cat-GV}
  Let \((\cat{C}, d)\) be a hom-finite left \GV{} category with dualising functor \(D_{\ell}\).
  Then \(\Cpsh_{\fin}\) can be endowed with the structure of a right \GV{} category  with dualising object
  \(\rd*{\cat{C}(\blank, d)} \in \Cpsh_{\fin}\) and dualising functor
  \begin{equation} \label{eq:dual-fun-cat}
    \mathsf{D}_r \from \Cpsh_{\fin}^{\op} \to \Cpsh_{\fin}, \qquad \qquad
    F \mapsto {F(D_{\ell}\blank)}^{*}.
  \end{equation}
\end{proposition}
\begin{proof}
  Using Equation~\eqref{eq:gv-dualising-from-internal-hom-vice-versa}, one can recover the right dualising functor of a \GV{} category from its right internal hom
  by evaluating it at the (potential) dualising object;
  \ie, \(\mathsf{D}_r \cong {[\blank, {\cat{C}(\blank, d)}^{*}]}_r\).
  In our case, we have
  \begin{align*}
    {[F, \rd*{\cat{C}(\blank, d)}]}_r x
    \overset{\eqref{eq:right-internal hom-simplification}}&{\cong} \!\!\int_b \kvect\big(F{[x,b]}_{\ell}, \kvect(\cat{C}(b,d), \k)\big) \cong \!\! \int_b \kvect\big(\cat{C}(b, d) \otimes F{[x,b]}_{\ell}, \k\big) \\
                                                          &\cong \kvect\Big(\!\int^b\!\!\cat{C}(b,d)\otimes F{[x,b]}_{\ell}, \k \Big)
                                                            \overset{\eqref{eq:coYoneda-covariant}}{\cong}  \kvect\big(F{[x,d]}_{\ell}, \k\big)
                                                            \overset{\eqref{eq:gv-dualising-from-internal-hom-vice-versa}}{\cong} \rd*{F(D_{\ell}x)}.
  \end{align*}
  To conclude the proof, notice that  \(\mathsf{D}_r \from \Cpsh_{\fin}^{\op}\to \Cpsh_{\fin}\) has a quasi-inverse
  \begin{equation*}
    \mathsf{D}_r^{-1} \from\Cpsh_{\fin} \to \Cpsh_{\fin}^{\op}, \qquad \qquad F \mapsto \rd*{F(D_{\ell}^{-1}\blank)}.\qedhere
  \end{equation*}
\end{proof}

While we will develop efficient means to detect tensor-representability and rigidity based on properties which exist in the abelian case, it is worthwhile to explore how such structures arise as an interplay between the duality type of the `base category' and that of \(\kvect\).

\begin{corollary}\label{cor:day-compact}
  Let \((\cat{C}, d)\) be a hom-finite left \GV{} category
  with dualising functor \(D_{\ell}\).
  The category \(\Cpsh_{\fin}\) is right tensor-representable if there are isomorphisms
  \begin{equation} \label{eq:trep-nice-enough}
    D_{\ell}^2a \cong a,
    \qquad \qquad D_{\ell}(a\otimes x) \cong D_{\ell}x \otimes D_{\ell}a,
    \qquad \qquad \text{ natural in } a,x \in \cat{C},
  \end{equation}
  and for all objects \(F, G \in \Cpsh_{\fin}\) we have
  \begin{equation} \label{eq:end-coend-commute}
    \int^a \kvect(F{[\blank,a]}_{\ell}, Ga) \cong \int_a \kvect(F{[\blank,a]}_{\ell},Ga).
  \end{equation}
\end{corollary}
\begin{proof}
  Let \(\mathsf{D}_r\) be the dualising functor of Proposition~\ref{prop:fun-cat-GV}.
  Fix objects \(F, G \in \Cpsh_{\fin}\) and compute
  \begin{align*}
    (\mathsf{D}_r F \ostar G) x
    \overset{\eqref{eq:day-convolution-simp-left}}%
    &{\cong}
      \int^a\!\!\mathsf{D}_r F({[a,x]}_{\ell}) \otimes Ga
      \overset{\eqref{eq:gv->internal hom}}{\cong}
      \int^a\!\!\mathsf{D}_r F(D_{\ell}(a \otimes D_{\ell}x)) \otimes Ga \\
    \overset{\ref{prop:fun-cat-GV}}%
    &{\cong}
      \int^a\!\!\rd*{F(D_{\ell}^2(a \otimes D_{\ell}x))} \otimes Ga
      \overset{\eqref{eq:trep-nice-enough}}{\cong}
      \int^a\!\!\kvect\big(F(a \otimes D_{\ell}x), Ga\big) \\
    \overset{\eqref{eq:trep-nice-enough}}%
    &{\cong}
      \int^a\!\!\kvect\big(F(D^2_{\ell} a \otimes D_{\ell} x), Ga\big)
      \overset{\eqref{eq:trep-nice-enough}}{\cong}
      \int^a\!\!\kvect\big(F(D_{\ell}(x \otimes D_{\ell} a)), Ga\big) \\
    \overset{\eqref{eq:gv->internal hom}\,+\,\eqref{eq:trep-nice-enough}}%
    &{\cong}
      \int^a\!\!\kvect\big(F({[x,a]}_{\ell}), Ga\big)
      \overset{\eqref{eq:end-coend-commute}}{\cong}
      \int_a\kvect\big(F({[x,a]}_{\ell}), Ga\big)
      \overset{\eqref{eq:left-internal hom-simplification}}
      = {[F,G]}_r x.
  \end{align*}
\end{proof}

\begin{remark}\label{rmk:conditions-for-early-properties}
  Before we discuss conditions for the interchangeability of ends and coends, let us briefly mention some cases where the dualising functor of a left \GV{} category \((\cat{C},d)\)
  admits natural isomorphisms as stated in Equation~\eqref{eq:trep-nice-enough}.

  The requirement \(D_{\ell}(x \otimes y) \cong D_{\ell}y \otimes D_{\ell}x\) implies that \(\cat{C}\) must be left tensor-repre\-sentable,
  since then for all \(x,y, z \in\cat{C}\) we have
  \[
    \cat{C}(x \otimes y, z)
    \overset{\eqref{eq:left-gv-condition-and-inverse}}{\cong}
    \cat{C}(x \otimes y \otimes D_{\ell}^{-1}z, d)
    \overset{\eqref{eq:left-gv-condition-and-inverse}}{\cong}
    \cat{C}(x, D_{\ell}(y \otimes D_{\ell}^{-1}z))
    \overset{\eqref{eq:trep-nice-enough}}{\cong}
    \cat{C}(x,z \otimes  D_{\ell}y).
  \]
  Notice, however, that in the absence of further coherence assumptions, this does not imply that \(\cat{C}\) must be rigid.

  The condition \(D_{\ell}^2 \cong \Id_{\cat{C}}\) is met for example if \(\cat{C}\) is braided since then for all \(x, y \in\cat{C}\) we have
  \[
    \cat{C}(x, y)
    \overset{\eqref{eq:left-gv-condition-and-inverse}}{\cong}
    \cat{C}(D_{\ell}y \otimes x,d)
    \cong
    \cat{C}(x \otimes D_{\ell}y  ,d)
    \overset{\eqref{eq:left-gv-condition-and-inverse}}{\cong}
    \cat{C}(x, D_{\ell}^2 y),
  \]
  and the Yoneda lemma  implies that \(\Id_{\cat{C}} \cong D_{\ell}^2\).
\end{remark}

The following constructions are an adaptation of~\cite{day06:compac}.
For simplicity, we replaced the promonoidal structure in \emph{loc.~cit.}
with the one induced from a monoidal structure.\footnote{\,%
  That is,
  if \(J\) and \(P\) are the promonoidal unit and multiplication,
  we set
  \begin{equation*}
    J \eqdef \cat{C}(1, \blank ) \qquad\qquad\text{and}\qquad\qquad P(x, y, z) \eqdef \cat{C}(x \otimes y, z).
  \end{equation*}%
}
\begin{remark}\label{rmk:preconditions-day-yoneda}
  Let \(\cat{C}\) be a category.
  Given any functor \(T\from \cat{C}^{\op} \otimes_{\k} \cat{C} \to \kVect\), there exists a canonical map
  \begin{equation}\label{eq:canonical-coend-map}
    \int^{a \in \cat{C}}\!\!\!\!\int_{b \in \cat{C}} \cat{C}(b, a) \otimes T(a, b)
    \to \int_{b \in \cat{C}} \int^{a \in \cat{C}} \cat{C}(b, a) \otimes T(a, b).
  \end{equation}
  This follows from the fact that ends and coends are functorial;
  in particular, for any arrow \(f \from x \to y\), the following diagram commutes:
  % file:///home/slot/repos/quiver/src/index.html?q=WzAsNSxbMCwwLCJcXGludF95IFxcbWF0aGNhbHtDfSh5LCBhKSBcXGNkb3QgVChiLCB5KSJdLFsxLDAsIlxcaW50X3kgXFxtYXRoY2Fse0N9KHksIGIpIFxcY2RvdCBUKGIsIHkpIl0sWzAsMSwiXFxpbnRfeSBcXG1hdGhjYWx7Q30oeSwgYSkgXFxjZG90IFQoYSwgeSkiXSxbMSwxLCJcXGludF54IFxcaW50X3kgXFxtYXRoY2Fse0N9KHksIHgpIFxcY2RvdCBUKHgsIHkpIl0sWzIsMiwiXFxpbnRfeSBcXGludF54IFxcbWF0aGNhbHtDfSh5LCB4KSBcXGNkb3QgVCh4LCB5KSJdLFswLDEsIlxcaW50X3kgXFxtYXRoY2Fse0N9KHksIGYpIFxcY2RvdCBUKGEsIHkpIl0sWzAsMiwiXFxpbnRfeSBcXG1hdGhjYWx7Q30oeSwgYSkgXFxjZG90IFQoZiwgeSkiLDJdLFsxLDMsIlxcYWxwaGEiXSxbMiwzLCJcXGFscGhhIiwyXSxbMyw0LCIiLDAseyJzdHlsZSI6eyJib2R5Ijp7Im5hbWUiOiJkb3R0ZWQifX19XSxbMiw0LCJcXGludF95IFxcYWxwaGEiLDIseyJjdXJ2ZSI6M31dLFsxLDQsIlxcaW50X3kgXFxhbHBoYSIsMCx7ImN1cnZlIjotM31dXQ==
  \begin{equation*}
    \begin{tikzcd}[column sep=large]
      {\int_b \cat{C}(b, x) \otimes T(y, b)} &[10pt] {\int_b \cat{C}(b, y) \otimes T(y, b)} \\
      {\int_b \cat{C}(b, x) \otimes T(x, b)} &[10pt] {\int^a\!\! \int_b \cat{C}(b, a) \otimes T(a, b)} \\
      &&[-30pt] {\int_b \int^a \cat{C}(b, a) \otimes T(a, b)}
      \arrow["{\int_b \cat{C}(b, f) \otimes T(y, b)}", from=1-1, to=1-2]
      \arrow["{\int_b \cat{C}(b, x) \otimes T(f, b)}"', from=1-1, to=2-1]
      \arrow["\can", from=1-2, to=2-2]
      \arrow["\can"', from=2-1, to=2-2]
      \arrow["{\exists!}", dotted, from=2-2, to=3-3]
      \arrow["{\int_b \can}"', curve={height=18pt}, from=2-1, to=3-3]
      \arrow["{\int_b \can}", curve={height=-18pt}, from=1-2, to=3-3]
    \end{tikzcd}
  \end{equation*}
\end{remark}

\begin{lemma}[{\cite[p.~1]{day06:compac}}]\label{lem:day-yoneda}
  Let \(\cat{C}\) be a hom-finite category,
  \(T \from \cat{C}^{\op} \otimes_{\k} \cat{C} \to \kvect\) a functor,
  and suppose that the map given in  Equation~\eqref{eq:canonical-coend-map} is invertible.
  If there is a natural isomorphism  \(\varphi \from \cat{C}(b, a) \iso \rd*{\cat{C}(a, b)}\),
  then
  \begin{equation*}
    \int^{a} T(a, a) \iso \int_{b} T(b, b).
  \end{equation*}
\end{lemma}
\begin{proof}
  We calculate
  \begin{align*}
    \int^{a} T(a, a)
    & \cong \int^{a}\!\!\!\int_{b} \kvect(\cat{C}(a, b), T(a, b))
      \cong \int^{a}\!\!\!\int_{b}  \rd*{\cat{C}(a, b)}\otimes T(a, b) \\
    \overset{\varphi^{-1}}&{\cong} \int^{a}\!\!\!\int_{b}  \cat{C}(b, a)\otimes T(a, b)
                            \overset{\eqref{eq:canonical-coend-map}}{\cong}  \int_{b} \int^{a} \cat{C}(b, a) \otimes T(a, b)
                            \cong \int_{b} T(b, b).
  \end{align*}
  Here, the first and last isomorphisms follow from Lemma~\ref{lem:enriched-yoneda} and the second one is a consequence of \(\cat{C}(a, b)\) being finite-dimensional and \(\kvect\) being rigid.
\end{proof}

In general, it seems difficult to decide whether there exists a natural isomorphism \(\varphi \from \cat{C}(b,a) \to \rd*{\cat{C}(a,b)}\).
However, for certain classes of examples, trace maps provide us with viable candidates.
\begin{remark}\label{rem:day-morphism-invertible}
  Let  \(\cat{C}\) be hom-finite and \emph{pivotal}—%
  that is, rigid monoidal such that there is a natural monoidal isomorphism \(\psi\from \rd{(\blank)}\to \ld{(\blank)}\) between the left and right dualising functors.
  If moreover \(\cat{C}(1,1) \cong \k\),
  then for all \(a, b \in \cat{C}\) we can consider:
  \begin{equation}\label{eq:day-canonical-candidate}
    \begin{gathered}
      \varphi_{a,b} \from
      \cat{C}(b,a) \to \rd*{\cat{C}(a,b)}, \qquad
      \varphi_{a,b}(f)g = \tr(fg)\eqdef \ev^{\ell}_{a}(\psi_{a}\otimes fg)\coev^{r}_{a}.
    \end{gathered}
  \end{equation}
  Its dual is
  \[
    \rd*{(\varphi_{a,b})}\from \cat{C}(a,b)\cong \rrd*{\cat{C}(a,b)}\to \rd*{\cat{C}(b,a)}, \quad
    \rd*{(\varphi_{a,b})}(g)f=\tr(fg)=\tr(gf)=\varphi_{b,a}(g)f.
  \]
  Thus, \(\varphi_{a,b}\) is injective if and only if \(\varphi_{b,a}\) is surjective.

  Suppose \(\k\) has characteristic zero and  \(\cat{C} \eqdef \lmod{H}\) is the category of finite-dimensional modules of a semisimple Hopf algebra \(H\).
  As a consequence of~\cite[Theorem~4]{larson-radford1988:SemisimpleCosemisimpleHopf}, left and right duals  coincide and  we can choose the `quantum trace' of Equation~\eqref{eq:day-canonical-candidate} to  agree with the usual trace of endomorphisms between finite-dimensional vector spaces.
  Let \(f\from M \to N\) be a morphism in \(\cat{C}\).
  Due to semisimplicity, the following short exact sequence splits.
  \begin{equation*}
    \begin{tikzcd}[ampersand replacement=\&]
      0 \&\& {\ker f} \&\& M \&\& {\im f} \&\& 0
      \arrow[from=1-1, to=1-3]
      \arrow[from=1-3, to=1-5]
      \arrow["f", from=1-5, to=1-7]
      \arrow["\iota"{description}, curve={height=-18pt}, dashed, from=1-7, to=1-5]
      \arrow[, from=1-7, to=1-9]
    \end{tikzcd}
  \end{equation*}
  The  arrow \(g= \iota \oplus 0 \from N \cong \im f \oplus N/\im f \to M\) satisfies
  \begin{equation*}
    \tr(fg) = \tr(f\iota)= \tr(\id_{\im f})=\dim \im f,
  \end{equation*}
  implying that \(\varphi_{M,N}\) is injective.

  Let us now consider a field  \(\k\) of characteristic \(p\) and \(q\in \mathbb{N}\) a multiple of \(p\).
  The group algebra \(H=\k[\GL_q(p)]\) of the group of invertible \(q\times q\)-matrices over \(\k\) provides us with an example of a pivotal hom-finite category where the morphism of Equation~\eqref{eq:day-canonical-candidate} is not invertible.
  Matrix-vector multiplication turns \(\k^q\) into a simple \(H\)-module whose endomorphism algebra is one-dimensional.
  This implies  for all \(f,g \in \End_{\GL_{q}(p)}(\k^p)\) that there is a \(\lambda \in \k\) such that \(\tr(fg)=\lambda \tr(\id_{\k^q}) = \lambda \cdot 0 = 0\).
\end{remark}

\begin{remark}\label{rmk:reduction-of-assumptions}
  In practise, one can use a combination of Corollary~\ref{cor:day-compact} and Lemma~\ref{lem:day-yoneda} to deduce that a given category of finite-dimensional functors is right tensor-representable.
  To reduce the number of axioms to be checked, note that the square of the left dualising functor being isomorphic to the identity is implied by the other assumptions.
  Let \((\cat{C},d)\) be a hom-finite left \GV{} category with dualising functor \(D_{\ell}\),
  such that there are isomorphisms
  \begin{equation*}
    \cat{C}(a,b) \cong \rd*{\cat{C}(b, a)}, \qquad \qquad
    D_{\ell}(a \otimes b) \cong D_{\ell}b \otimes D_{\ell}a ,\qquad \qquad \text{ natural in } a, b \in\cat{C}.
  \end{equation*}
  Using that
  \begin{align*}
    \cat{C}(d,a \otimes D_{\ell}b)
    \overset{\eqref{eq:left-gv-condition-and-inverse}}%
    &{\cong} \cat{C}(d \otimes D_{\ell}^{-1} (a \otimes D_{\ell} b), d)
      \cong \cat{C}(d \otimes b \otimes D_{\ell}^{-1}a, d) \\
    \overset{\eqref{eq:left-gv-condition-and-inverse}}%
    &{\cong} \cat{C}(D_{\ell}^{-1}a, D_{\ell}^{-1}(d \otimes b))
      \overset{\eqref{eq:gv:1:from:d:and:d:from:1}}%
      {\cong} \cat{C}(D_{\ell}^{-1}a, D_{\ell}^{-1}b)
      \cong \cat{C}(b,a),
  \end{align*}
  one obtains
  \begin{equation*}
    \cat{C}(a, D_{\ell}^2b)
    \overset{\eqref{eq:left-gv-condition-and-inverse}}%
    {\cong} \cat{C}(a \otimes D_{\ell}b, d)
    \cong \rd*{\cat{C}(d, a \otimes D_{\ell}b)}
    \cong \rd*{\cat{C}(b, a)}
    \cong \rrd*{\cat{C}(a,b)}
    \cong \cat{C}(a,b).
  \end{equation*}
  The claim follows from the Yoneda lemma.
\end{remark}

\begin{lemma}\label{lem:semisimple-lets-ends-commute}
  Assume \(\k\) to be a perfect field.
  If the hom-finite category \(\cat{C}\) has finitely many objects
  and \(\Cpsh_{\fin}\)  is semisimple,
  then the canonical map
  \begin{equation*}
    \int^a \!\!\!\!\int_b T(b,b,a,a) \cong \int_b \int^a T(a,a,b,b)
  \end{equation*}
  is invertible for all functors \(T \from \cat{C}^{\op}\otimes_{\k}\cat{C} \otimes_{\k}\cat{C}^{\op}\otimes_{\k}\cat{C} \to \kvect\).
\end{lemma}
\begin{proof}
  We endow the vector space \(A \eqdef \bigoplus_{a, b\in\cat{C}}\cat{C}(a,b)\) with the structure of an associative unital algebra via the multiplication, specified on elements \(f \in \cat{C}(a, b),\ g \in \cat{C}(c, d)\) by
  \begin{equation*}
    f\cdot g =
    \begin{cases*}
      g f, & for \(b = c\), \\
      0, & otherwise.
    \end{cases*}
  \end{equation*}
  The unit \(\sum_{a \in\cat{C}} \id_a\) of \(A\)  is a sum of orthogonal idempotents.
  Thus, any right module \(M\) of \(A\) decomposes as a vector space into a direct sum \(M \cong \bigoplus_{a \in\cat{C}} M_a\), where \(M_a = M\cdot\id_a\), and the action of any  \(f \in\cat{C}(a,b)\)  defines a linear map \(M_a \to M_b\).
  Accordingly, a morphism of right \(A\)-modules corresponds to a collection of homomorphisms \({\{\,\phi_a \from M_a \to N_a\,\}}_{a \in\cat{C}}\) such that
  \begin{equation*}
    \phi_b(m \cdot f) = \phi_a(m) \cdot f,
    \qquad \qquad \text{ for all }
    m \in M_a, f\in \cat{C}(a,b).
  \end{equation*}
  This defines a \(\k\)-linear functor \(\Theta \from \rMod{A} \to \Cpsh\).
  Its quasi-inverse \(\Omega \from \Cpsh \to \rMod{A}\) maps any \(F \in \Cpsh\) to the module \(\bigoplus_{a \in\cat{C}}Fa\), and any  arrow \({\{\,\psi_a \from Fa \to Ga\,\}}_{a \in\cat{C}}\) to the module homomorphism \(\bigoplus_{a\in\cat{C}}\psi_a \from \bigoplus_{a \in\cat{C}} Fa \to \bigoplus_{a\in\cat{C}} Ga\).
  Via this identification \(\Cpsh_{\fin}\) corresponds to the category \(\rmod{A}\) of finite-dimensional right \(A\)-modules.

  Since \(\Cpsh_{\fin}\) is semisimple, so is \(A\) and \(A^{\op}\).
  Furthermore, \(\k\) being perfect implies that the tensor product of any two finite-dimensional semisimple algebras is semisimple, see~\cite[Section~3]{farb-dennis1993:Noncommutative}.
  In particular, \(B \eqdef A^{\op}\otimes_{\k}A \otimes_{\k}A^{\op}\otimes_{\k}A\) is semisimple
  and a similar argument as before shows that \([\cat{C}^{\op}\otimes_{\k}\cat{C} \otimes_{\k}\cat{C}^{\op}\otimes_{\k}\cat{C}, \kvect] \cong \rmod{B}\).
  We write \(\cat{D} \eqdef [\cat{C}^{\op}\otimes_{\k} \cat{C}, [\cat{C}^{\op} \otimes_{\k}\cat{C}, \kvect]]\).

  There is a \(\k\)-linear equivalence of categories
  \begin{gather*}
    [\cat{C}^{\op}\otimes_{\k} \cat{C} \otimes_{\k}\cat{C}^{\op} \otimes_{\k}\cat{C}, \kvect] \to \cat{D},\\
    F \mapsto \widehat F, \qquad \text{where} \qquad \widehat F(x,y)[u,v] = F(u,v,x,y).
  \end{gather*}
  The functor \(\text{end} \from\cat{D}\to [\cat{C}^{\op} \otimes_{\k}\cat{C}, \kvect]\), given by \(\text{end}(F)(u, v) = \left(\int_b F(b,b)\right)(u, v)\),
  is left exact since limits commute with limits.
  Due to semisimplicity, any short exact sequence in \(\cat{D}\) splits and split epimorphisms are preserved by all functors.
  Thus, the functor \(\text{end} \from\cat{D}\to [\cat{C}^{\op} \otimes_{\k}\cat{C}, \kvect]\) is exact and must preserve colimits.

  Given \(F\in  [\cat{C}^{\op}\otimes_{\k}\cat{C} \otimes_{\k}\cat{C}^{\op}\otimes_{\k}\cat{C}, \kvect]\), we now compute
  \begin{equation*}
    \int^a\!\!\!\!\int_b F(a,a)(b,b)
    \cong \int^a\!\!\!\!\int_b\widehat{F}(b,b)[a,a]
    \cong \int_b\int^a\widehat{F}(b,b)(a,a)
    \cong \int_b\int^a F(a,a,b,b).
  \end{equation*}
\end{proof}

\subsection{Tensor-representability, rigidity, and finitely-generated projective functors}\label{sec: trep-rigid-fgp}
Let us conclude this section by stating explicit criteria for certain functor categories to be tensor-representable or even rigid.
As with modules over (commutative) rings both notions are closely connected with objects being finitely-generated projective.
Again, we implicitly assume all categories and functors to be \(\k\)-linear, for a field \(\k\).
However, we do not require  finiteness of hom-spaces and work solely under the assumptions made in Hypothesis~\ref{hypo:eso-enriched}.

The following lemma is a special case of conditions discussed in~\cite[Chapter~10]{prest09:purit}.

\begin{lemma}\label{lem:projective-if-representable}
  Let \(\cat{C}\) be a  category.
  For any functor \(F \in \Cpsh\)  the following are equivalent:
  \begin{thmlist*}
    \item \(F\) is finitely-generated projective,\label{itm:fgp-char-fgp}
    \item \(F\) is a direct summand of a finite direct sum of representable functors, and\label{itm:fgp-char-dir-sum}
    \item the functor \(\Cpsh(F, \blank)\) commutes with small colimits.\label{itm:fgp-char-colim}
  \end{thmlist*}
\end{lemma}
\begin{proof}
  The equivalence between~\ref{itm:fgp-char-fgp} and~\ref{itm:fgp-char-dir-sum} is proven in~\cite[Corollary~10.1.14]{prest09:purit}.
  In order to show that~\ref{itm:fgp-char-dir-sum} and~\ref{itm:fgp-char-colim} are equivalent, observe that the full subcategory \(\overline{\cat{C}}\) of \(\Cpsh\) consisting of direct summands of finite direct sums of representable functors is Cauchy complete\footnote{\,%
    A category \( \cat{D}\) is \emph{Cauchy complete} if it has all \emph{absolute colimits}. These are colimits which are preserved by all functors with source category \( \cat{D}\).
  } by~\cite[Corollary~4.22]{lack-tendas2022:FlatFilteredColimts}.
  The claim now follows by proceeding analogous to Proposition~2 of~\cite{borceux-dejean1986:CauchyCompletionCategoryTheory}.
\end{proof}

For a category \(\cat{C}\),
let us write \(\overline{\cat{C}}\) for the full subcategory of \(\Cpsh\)
whose objects are direct summands of finite direct sums of representable functors.
In accordance to~\cite[Section~5.5]{kelly05:basic}, we refer to it as the \emph{Cauchy completion} of \(\cat{C}^{\op}\).
In order to show that \(\overline{\cat{C}}\) is monoidal and investigate the types of dualities induced by \(\cat{C}^{\op}\), we need a compatibility between the convolution tensor product and the Yoneda embedding due to Day,~\cite{day71:const},
see also~\cite{im86:univ-conv}.

\begin{lemma}\label{lem:yoneda-monoidal}
  For a  monoidal category \(\cat{C}\),
  the Yoneda embedding \(\yo \from \cat{C}^{\op} \to \Cpsh\) is a strong monoidal functor.
\end{lemma}

In view of Lemma~\ref{lem:projective-if-representable}, we shall adopt the following notation
in order to emphasise that a given functor \(X \in \overline{\cat{C}}\) is a direct summand of a direct sum of representables:
\[
  \begin{tikzcd}[cramped]
    X & {\bigoplus\limits_{i=1}^n \yo(u_i),}
    \arrow["{\iota_X}", shift left, from=1-1, to=1-2]
    \arrow["{\pi_X}", shift left, from=1-2, to=1-1]
  \end{tikzcd}
\]
where \(\iota_X\) and \(\pi_X\) are sections and retracts of each other; \ie, \(\pi_X \circ \iota_X = \id_X\).

We may now characterise the closed monoidal structure of the Cauchy completion.

\begin{proposition}\label{prop:cauchy-completion-is-trep}
  The Cauchy completion \(\overline{\cat{C}}\) of a  right closed monoidal category \(\cat{C}^{\op}\) is right closed monoidal.
  For all objects \(
  \begin{tikzcd}[cramped]
    X & {\bigoplus_{i=1}^n \yo(u_i)}
    \arrow["{\iota_X}", shift left, from=1-1, to=1-2]
    \arrow["{\pi_X}", shift left, from=1-2, to=1-1]
  \end{tikzcd}
  \) and
  \(
  \begin{tikzcd}[cramped]
    Y & {\bigoplus_{j=1}^m \yo(v_i)}
    \arrow["{\iota_Y}", shift left, from=1-1, to=1-2]
    \arrow["{\pi_Y}", shift left, from=1-2, to=1-1]
  \end{tikzcd}
  \),
  the right internal hom satisfies
  \begin{equation}\label{eq:char-of-internal hom-of-trep}
    \begin{tikzcd}[ampersand replacement=\&]
      {[X,Y]}_r \& \& {[\oplus_{i=1}^n \yo u_i, \oplus_{j=1}^m \yo v_j]\cong \bigoplus\limits_{i=1}^n \bigoplus\limits_{j=1}^m \yo({[u_i,v_j]}_r).}
      \arrow["{{[\pi_{X},\iota_Y]}_r}", shift left, from=1-1, to=1-3]
      \arrow["{{[\iota_X, \pi_Y]}_r}", shift left, from=1-3, to=1-1]
    \end{tikzcd}
  \end{equation}
\end{proposition}
\begin{proof}
  The fact that \(\overline{\cat{C}} \subseteq \Cpsh\) is closed under taking tensor products is due to the Yoneda embedding being strong monoidal, and \(\overline{\cat{C}}\) being the full subcategory of \(\Cpsh\) whose objects are direct summands of finite direct sums of representables.

  Let us now compute the right internal hom \({[X,Y]}_r \in \Cpsh\) of two objects \(X, Y \in \overline{\cat{C}}\).
  We may write \(X= \colim_i\cat{C}(u_i, \blank)\) and \(Y = \colim_j\cat{C}(v_j, \blank)\) as direct summands of finite direct sums of representables.
  A straightforward computation yields
  \begin{align*}
    {[X,Y]}_r
    & \overset{~\eqref{eq:Day-internal-hom-right}}{=} \int_{a, b}\kVect(\cat{C}(a \otimes \blank, b), \kVect(Xa, Yb))
      \overset{\eqref{eq:coYoneda-covariant}}{\cong} \int_{a}\kVect(Xa, Y(a\otimes \blank))\\
    & \cong \int_{a} \kVect(\colim_i\cat{C}(u_i, a), \colim_j \cat{C}(v_j, a \otimes \blank)) \\
    & \cong \lim\nolimits_i \colim_j  \int_{a}\kVect(\cat{C}(u_i, a), \cat{C}(v_j, a\otimes \blank))
      \overset{\eqref{eq:coYoneda-covariant}}\cong \lim\nolimits_i \colim_j\cat{C}(v_j, u_i \otimes \blank) \\
    & \cong \lim\nolimits_i \colim_j\cat{C}({[u_i, v_j]}_r, \blank).
  \end{align*}
  It follows that \(\overline{\cat{C}}\) is right closed monoidal and that Equation~\eqref{eq:char-of-internal hom-of-trep} holds.
\end{proof}

The previous result can be understood as a variation of the fact that a direct summand or direct sum of (rigidly) dualisable objects is dualisable.
Indeed, suppose \(\cat{C}^{\op}\) to be a  closed monoidal category
and consider three direct summands
\adjustbox{trim=0em 0.5em 0em 0em}{\(\begin{tikzcd}[ampersand replacement=\&, cramped]
    X \& {U}
    \arrow["{\iota_X}", shift left, from=1-1, to=1-2]
    \arrow["{\pi_X}", shift left, from=1-2, to=1-1]
  \end{tikzcd}
  \)},
\adjustbox{trim=0em 0.5em 0em 0em}{\(
  \begin{tikzcd}[ampersand replacement=\&, cramped]
    Y \& V
    \arrow["{\iota_Y}", shift left, from=1-1, to=1-2]
    \arrow["{\pi_Y}", shift left, from=1-2, to=1-1]
  \end{tikzcd}
  \)},
\adjustbox{trim=0em 0.5em 0em 0.35em}{\(\begin{tikzcd}[ampersand replacement=\&, cramped]
    Z \& {W}
    \arrow["{\iota_Z}", shift left, from=1-1, to=1-2]
    \arrow["{\pi_Z}", shift left, from=1-2, to=1-1]
  \end{tikzcd}\)} of objects in \(\overline{\cat{C}}\).
The following diagram, whose horizontal arrows are the isomorphisms of the tensor-hom adjunction of \( \overline{\cat{C}}\), commutes:
\begin{equation} \label{eq:comm-diag-for-unit-and-counit-summand}
  \begin{tikzcd}[ampersand replacement=\&]
    {\overline{\cat{C}}(X \ostar Y, Z)} \&[3em]\& {\overline{\cat{C}}(Y,{[X,Z]}_{r})} \\[1em]
    {\overline{\cat{C}}(U \ostar V, W)} \&\& {\overline{\cat{C}}(V, {[U,W]}_{r})}
    \arrow["{\overline{\cat{C}}(\pi_X \ostar \pi_Y, \iota_Z)}"', shift right, from=1-1, to=2-1]
    \arrow["{\phi_{U,V,W}}", from=2-1, to=2-3]
    \arrow["{\phi_{X,Y,Z}}", from=1-1, to=1-3]
    \arrow["{\overline{\cat{C}}(\pi_V,{[\pi_X,\iota_Z]}_r)}"', shift right, from=1-3, to=2-3]
    \arrow["{\overline{\cat{C}}(\iota_X \ostar \iota_Y, \pi_Z)}"', shift right, from=2-1, to=1-1]
    \arrow["{\overline{\cat{C}}(\iota_V,{[\iota_X,\pi_Z]}_r)}"', shift right, from=2-3, to=1-3]
  \end{tikzcd}
\end{equation}

Thus the unit and counit of the adjunction \(\adj{X \ostar \blank}{{[X, \blank]}_r}{\overline{\cat{C}}}{\overline{\cat{C}}}\) satisfy
\begin{equation} \label{eq:unit-counit-formula-direct-summand}
  \eta^{(X)}_Y = [\iota_X, \pi_X \ostar \pi_Y] \eta^{(U)}_V \iota_Y \qquad\text{ and }\qquad
  \varepsilon^{(X)}_Y = \pi_Y\varepsilon^{(U)}_V(\iota_X\ostar [\pi_X,\iota_Y]).
\end{equation}

We will now prove that \(\cat{C}^{\op}\) and its Cauchy completion share the same type of duality.

\begin{corollary}\label{cor:pres-ref-dual-struct}
  Let \(\cat{C}^{\op}\) be a right closed monoidal category.
  We have:
  \begin{thmlist*}
    \item \(\cat{C}^{\op}\) is right rigid if and only if \(\overline{\cat{C}}\) is right rigid.\label{itm:rigid-implies-rigid}
    \item \(\cat{C}^{\op}\) is right tensor-representable if and only if \(\overline{\cat{C}}\) is right tensor-representable.\label{itm:trep-implies-trep}
    \item \((\cat{C}^{\op},d)\) is a right \GV{} category if and only if \((\overline{\cat{C}}, \yo(d))\) is.\label{itm:GV-implies-GV}
  \end{thmlist*}
\end{corollary}

\begin{proof}
  That right rigidity and tensor-representability of \(\cat{C}^{\op}\) imply the same property for \(\overline{\cat{C}}\) follows straightforwardly from the description of the internal hom of \( \overline{\cat{C}}\) and the units and counits of the tensor-hom adjunctions, see Equation~\eqref{eq:unit-counit-formula-direct-summand}.

  Now suppose that \((\cat{C}^{\op}, d)\) is a right \GV{} category.
  Its right dualising functor is, up to natural isomorphism, given by \({[\blank, d]}_r \from\cat{C} \to\cat{C}^{\op}\).
  A direct computation using Proposition~\ref{prop:cauchy-completion-is-trep} shows that
  the functor \({[\blank, \yo(d)]}_r \from \overline{\cat{C}}^{\op}\to \overline{\cat{C}}\) is an equivalence,
  and therefore \(\big(\overline{\cat{C}}, \yo(d)\big)\) is a right \GV{} category.

  By Proposition~\ref{prop:cauchy-completion-is-trep} the right internal hom of two representable functors \(\yo(x)\) and \(\yo(y)\) is \({[\yo(x),\yo(y)]}_r \cong\yo({[x,y]}_r)\).
  Thus, the converse of any of the three statements is a consequence of the fact that, via the Yoneda embedding, \(\cat{C}^{\op}\) is equivalent as a right closed monoidal category to the full subcategory of \(\overline{\cat{C}}\) whose objects are representable functors.
\end{proof}

In many cases arising in representation theory,
like modules over commutative rings or \(\k\)-algebras,
rigidity and tensor-representability are equivalent;
see for example~\cite[Proposition~2.1]{niefield17:coexp} for a slightly more general statement.
Since a \(\k\)-algebra can be interpreted as a \(\k\)-linear category with one object,
the following proposition may be seen as a `many object version' of the classical case.

\begin{proposition}\label{prop:characterisation-of-rigidity}
  Let \(\cat{C}\) be a left rigid monoidal category and \(X \in \Cpsh\).
  The following properties for \(X\)  are equivalent:
  \begin{thmlist*}
    \item it has a right rigid dual,\label{itm:char-rigid}
    \item there exists a \(\mathbb{D}X \in \Cpsh\) such that \(\adj{X \ostar \blank}{\mathbb{D}X \ostar \blank}{\Cpsh}{\Cpsh}\), and\label{itm:char-trep}
    \item it is finitely-generated projective.\label{itm:char-fgp}
  \end{thmlist*}
\end{proposition}
\begin{proof}
  By definition,~\ref{itm:char-rigid} implies~\ref{itm:char-trep}.
  To show that~\ref{itm:char-trep} implies~\ref{itm:char-fgp}, we assume that there exists a \(\mathbb{D}X \in \Cpsh\)  such that  \(\adj{X \ostar \blank}{\mathbb{D}X \ostar \blank}{\Cpsh}{\Cpsh}\)  and fix a small colimit \(\colim_{i \in I} Fi \in \Cpsh\) of some diagram \(F \from I \to  \Cpsh\).
  For every \(i \in I\), write \(\iota_i \from Fi \to \colim_{i \in I} Fi \in \Cpsh\) for its structure morphisms.
  Now consider the commutative diagram:
  \[
    \begin{tikzcd}[ampersand replacement=\&,column sep=3.5em]
      {\colim_{i \in I} \Cpsh (X, Fi)} \&\& {\Cpsh (X, \colim_{i \in I} Fi)} \\
      \\
      {\colim_{i \in I} \Cpsh (1,\mathbb{D}X \ostar Fi)} \&\& { \Cpsh (1, \mathbb{D}X \ostar \colim_{i \in I} Fi)}
      \arrow["{\colim_{i \in I}(\Cpsh(X,\iota_i))}", from=1-1, to=1-3]
      \arrow["{\colim_{i \in I} \phi_{X, Fi}}", from=1-1, to=3-1]
      \arrow["{\phi_{X, \colim_{i \in I} Fi}}", from=1-3, to=3-3]
      \arrow["{\colim_{i \in I}(\Cpsh(1,\mathbb{D}X\ostar \iota_i))}", from=3-1, to=3-3]
    \end{tikzcd}
  \]
  Its horizontal arrows are induced by the universal property of the colimit and the vertical arrows are due to the tensor-hom adjunction of \(\Cpsh\).
  The functors \(\Cpsh(1, \blank)\) and \(\mathbb{D}X \ostar \blank\) commute with all small colimits;
  the first one due to the fact that \(\Cpsh(1, \blank)\) is finitely-generated projective, see Lemma~\ref{lem:projective-if-representable}, and the second one due to Day convolution being a colimit.
  Therefore, the horizontal arrow at the bottom is invertible.
  Since the vertical arrows are  also invertible, the canonical arrow \(\colim_{i \in I} \Cpsh (X, Fi) \to \Cpsh (X, \colim_{i \in I} Fi)\) displayed at the top of the diagram must be an isomorphism.
  Again, using Lemma~\ref{lem:projective-if-representable}, \(X\) must be finitely-generated projective.
  Thus,~\ref{itm:char-trep} implies~\ref{itm:char-fgp}.

  Finally, if \(X \in \Cpsh\) is finitely-generated projective, then it is contained in the Cauchy-completion of \(\cat{C}^{\op}\),
  which is right rigid.
  Then, by Corollary~\ref{cor:pres-ref-dual-struct}, so is \(\overline{\cat{C}}\).
  Therefore \(X\) admits a rigid dual and \(\ref{itm:char-fgp}\implies\ref{itm:char-rigid}\).
\end{proof}

\section{Applications in representation theory}\label{sec:mky-2-grps}
We conclude the article by applying the abstract machinery of the previous section to four kinds of representation theoretic examples.

First, we provide in Theorem~\ref{thm:sl2-is-trep} an example of a non-rigid tensor-representable \(r\)-category in terms of so-called \(\mathfrak{sl}_{2}(\mathbb{C})\)-crystals.

Then, we discuss the definition of Boolean algebras, some of their applications in group and ring theory, and  how they induce abelian, \(\k\)-linear \GV{} categories, see~Proposition~\ref{prop:qf2-algebras-from-boolean}.

Next, we focus on Mackey functors.
These are, roughly speaking, collections of vector spaces indexed by all subgroups of a fixed finite group together with morphisms subject to relations resembling the behaviour of `induction', `restriction', and `conjugation' operations, including the eponymous Mackey identity;
see~\cite{lindner76:mackey,thevenaz95:mackey}.
Finite-dimensional Mackey functors form a \GV{} category;
we show in Proposition~\ref{prop:mky-GV-and-rigid} that it is rigid if and only if it is semisimple.

The last example arises in the study of crossed modules,
which---in categorical terms---correspond to strict 2-groups.
The functors from any finite strict 2-group to \(\kvect\) form an abelian monoidal category
that is equivalent to a direct sum of representation categories of the isotropy-group of the monoidal unit of the 2-group.
The article is concluded with Proposition~\ref{prop:2-grp-GV-and-rigid}.
Therein, we prove the rigidity of this category to be equivalent to the semisimplicity of a certain group algebra.

\subsection{\texorpdfstring{\(\mathfrak{sl}_{2}\)}{sl2}-crystals}\label{sec:sl-2-crystals}
In the following, we work over the field \( \Bbbk=\mathbb{C}\).
As in~\cite{alqady-stroiński2025:TemperleyLieb}, we will define the \(\kVect[ \mathbb{C}]\)-enriched category of \(\mathfrak{sl}_{2}\)-crystals as the Cauchy completion of the so-called crystal Temperley--Lieb category.
The latter can be  constructed in terms of generators and relations.
For details of this type of construction we refer the reader to~\cite[Chapter~XII]{Kassel1998}.

\begin{definition}\label{def:tl-crys}
  The \emph{crystal Temperley--Lieb category} \(\mathsf{CTL}\) is the \(\kVect[\mathbb{C}]\)-enriched category generated by a single object \(t\) and the morphisms
  \begin{equation}\label{eq:generating-morphisms-tl}
    e: t\otimes t \to 1, \qquad\quad c \from 1 \to t\otimes t,
  \end{equation}
  subject to the relations
  \begin{subequations}
    \begin{gather}
      (e\otimes t)(t \otimes c) = 0 = (t \otimes e)(c \otimes t), \label{eq:tlc-snake}\\
      ec= 1. \label{eq:tlc-circle}
    \end{gather}
  \end{subequations}
\end{definition}

\begin{remark}\label{rmk:temperley--lieb-as-string-diags}
  The diagrammatics of the Temperley--Lieb category was first described in~\cite[Section~2]{khovanov1997:GraphicalKazhdanLusztig},
  see also~\cite{alqady-stroiński2025:TemperleyLieb} for the case of \(\mathfrak{sl}_2\)-crystals.
  Consider non-negative integers \(n, m \in \mathbb{N}_{0}\).
  The vector space \(\mathsf{CTL}(t^{\otimes m},t^{\otimes n})\) of morphisms between the \(m\)-fold and \(n\)-fold tensor power of \(t\) can be identified with the space of all linear combinations of unoriented (non-crossing) string diagrams \(S \subset \mathbb{R} \times [0,1]\) with
  \[
    S \cap \mathbb{R}\times \{0,1\} = \{1, \dots , m\}\times \{0\} \cup\{1, \dots , n\} \times\{1\}.
  \]
  A typical diagram consists of three types of strands: trough-strands, caps, and cups:
  \begin{center}
    \begin{tikzpicture}[scale =0.5]
      \begin{pgfonlayer}{nodelayer}
        \node [style=none] (0) at (3, 0) {};
        \node [style=none] (1) at (7, 0) {};
        \node [style=none] (2) at (1, 5) {};
        \node [style=none] (3) at (3, 5) {};
        \node [style=none] (4) at (12, 0) {};
        \node [style=none] (5) at (14, 0) {};
        \node [style=none] (6) at (14, 5) {};
        \node [style=none] (7) at (12, 5) {};
        \node [style=none] (8) at (10, 5) {};
        \node [style=none] (9) at (7, 5) {};
        \node [style=none] (10) at (5, 5) {};
        \node [style=none] (11) at (9, 0) {};
        \node [style=none] (12) at (5, 0) {};
        \node [style=none] (13) at (1, 0) {};
        \node [style=none] (14) at (8.5, 5) {\(\dots\)};
        \node [style=none] (15) at (10.5, 0) {\(\dots\)};
        \node [style=none] (17) at (0, 5.75) {};
        \node [style=none] (18) at (0.25, 0) {};
        \node [style=none] (19) at (-0.25, 0) {};
        \node [style=none] (20) at (0, -0.75) {};
        \node [style=none] (21) at (0.25, 5) {};
        \node [style=none] (22) at (-0.25, 5) {};
        \node [style=none] (23) at (-0.75, 0) {\(0\)};
        \node [style=none] (24) at (-0.75, 5) {};
        \node [style=none] (25) at (-0.75, 5) {\(1\)};
        \node [style=none] (26) at (14.75, -1) {};
        \node [style=none] (27) at (1, -0.75) {};
        \node [style=none] (28) at (1, -1.25) {};
        \node [style=none] (29) at (0.25, -1) {};
        \node [style=none] (30) at (3, -0.75) {};
        \node [style=none] (31) at (3, -1.25) {};
        \node [style=none] (32) at (1, -1.75) {\(1\)};
        \node [style=none] (33) at (3, -1.75) {};
        \node [style=none] (34) at (3, -1.75) {\(2\)};
        \node [style=none] (35) at (5, -0.75) {};
        \node [style=none] (36) at (5, -1.25) {};
        \node [style=none] (37) at (5, -1.75) {};
        \node [style=none] (38) at (5, -1.75) {\(3\)};
        \node [style=none] (39) at (7, -0.75) {};
        \node [style=none] (40) at (7, -1.25) {};
        \node [style=none] (41) at (7, -1.75) {};
        \node [style=none] (42) at (7, -1.75) {\(4\)};
        \node [style=none] (43) at (9, -0.75) {};
        \node [style=none] (44) at (9, -1.25) {};
        \node [style=none] (45) at (9, -1.75) {};
        \node [style=none] (46) at (9, -1.75) {\(5\)};
        \node [style=none] (47) at (12, -0.75) {};
        \node [style=none] (48) at (12, -1.25) {};
        \node [style=none] (49) at (12, -1.75) {};
        \node [style=none] (50) at (12, -1.75) {\(m-1\)};
        \node [style=none] (51) at (14, -0.75) {};
        \node [style=none] (52) at (14, -1.25) {};
        \node [style=none] (53) at (14, -1.75) {};
        \node [style=none] (54) at (14, -1.75) {\(m\)};
        \node [style=none] (55) at (14.75, 6) {};
        \node [style=none] (56) at (1, 5.75) {};
        \node [style=none] (57) at (1, 6.25) {};
        \node [style=none] (58) at (0.25, 6) {};
        \node [style=none] (59) at (3, 5.75) {};
        \node [style=none] (60) at (3, 6.25) {};
        \node [style=none] (61) at (1, 6.75) {\(1\)};
        \node [style=none] (62) at (3, 6.75) {};
        \node [style=none] (63) at (3, 6.75) {\(2\)};
        \node [style=none] (64) at (5, 5.75) {};
        \node [style=none] (65) at (5, 6.25) {};
        \node [style=none] (66) at (5, 6.75) {};
        \node [style=none] (67) at (5, 6.75) {\(3\)};
        \node [style=none] (68) at (7, 5.75) {};
        \node [style=none] (69) at (7, 6.25) {};
        \node [style=none] (70) at (7, 6.75) {};
        \node [style=none] (71) at (7, 6.75) {\(4\)};
        \node [style=none] (72) at (10, 5.75) {};
        \node [style=none] (73) at (10, 6.25) {};
        \node [style=none] (74) at (10, 6.75) {};
        \node [style=none] (75) at (10, 6.75) {\(n-2\)};
        \node [style=none] (76) at (12, 5.75) {};
        \node [style=none] (77) at (12, 6.25) {};
        \node [style=none] (78) at (12, 6.75) {};
        \node [style=none] (79) at (12, 6.75) {\(n-1\)};
        \node [style=none] (80) at (14, 5.75) {};
        \node [style=none] (81) at (14, 6.25) {};
        \node [style=none] (82) at (14, 6.75) {};
        \node [style=none] (83) at (14, 6.75) {\(n\)};
        \node [style=none] (84) at (15, 5.75) {};
        \node [style=none] (85) at (14.75, 0) {};
        \node [style=none] (86) at (15.25, 0) {};
        \node [style=none] (87) at (15, -0.75) {};
        \node [style=none] (88) at (14.75, 5) {};
        \node [style=none] (89) at (15.25, 5) {};
        \node [style=none] (90) at (15.75, 0) {\(0\)};
        \node [style=none] (91) at (15.75, 5) {};
        \node [style=none] (92) at (15.75, 5) {\(1\)};
        \node [style=none] (93) at (-2, 2) {a copy of \(e\)};
        \node [style=none] (94) at (18.75, 5.25) {a copy of \(c\)};
        \node [style=none] (95) at (-2, 1.5) {};
        \node [style=none] (96) at (2, 1.25) {};
        \node [style=none] (97) at (18.75, 4.75) {};
        \node [style=none] (98) at (13, 3.75) {};
        \node [style=none, align = left] (99) at (20.25, 1.75) {a through-strand is a \\ copy of the identity map};
        \node [style=none] (100) at (14, 1.25) {};
        \node [style=none] (101) at (20.25, 1) {};
      \end{pgfonlayer}
      \begin{pgfonlayer}{edgelayer}
        \draw [very thick, in=630, out=270, looseness=1.50] (6.center) to (7.center);
        \draw [very thick, in=270, out=90] (5.center) to (8.center);
        \draw [very thick, in=-90, out=-90, looseness=1.50] (10.center) to (9.center);
        \draw [very thick, in=90, out=-90] (3.center) to (4.center);
        \draw [very thick, in=90, out=90, looseness=1.50] (1.center) to (11.center);
        \draw [very thick, in=90, out=90, looseness=1.50] (13.center) to (0.center);
        \draw [very thick, in=270, out=90] (12.center) to (2.center);
        \draw (19.center) to (18.center);
        \draw (20.center) to (17.center);
        \draw (22.center) to (21.center);
        \draw (28.center) to (27.center);
        \draw (29.center) to (26.center);
        \draw (31.center) to (30.center);
        \draw (36.center) to (35.center);
        \draw (40.center) to (39.center);
        \draw (44.center) to (43.center);
        \draw (48.center) to (47.center);
        \draw (52.center) to (51.center);
        \draw (57.center) to (56.center);
        \draw (58.center) to (55.center);
        \draw (60.center) to (59.center);
        \draw (65.center) to (64.center);
        \draw (69.center) to (68.center);
        \draw (73.center) to (72.center);
        \draw (77.center) to (76.center);
        \draw (81.center) to (80.center);
        \draw (86.center) to (85.center);
        \draw (87.center) to (84.center);
        \draw (89.center) to (88.center);
        \draw [dashed, ->, in=90, out=-90, looseness=2.00] (95.center) to (96.center);
        \draw [dashed, ->, in=630, out=270] (97.center) to (98.center);
        \draw [dashed, ->, in=375, out=270] (101.center) to (100.center);
      \end{pgfonlayer}
    \end{tikzpicture}
  \end{center}
  In this diagrammatical language, placing diagrams next to each other corresponds to tensoring morphisms.
  Composition is realised by placing diagrams on top of each other and applying the defining relations.
  These take the form:
  \begin{center}
    \begin{tikzpicture}[scale=0.5]
      \begin{pgfonlayer}{nodelayer}
        \node [style=none] (2) at (1, 3) {};
        \node [style=none] (3) at (5, 5) {};
        \node [style=none] (7) at (18, 2.5) {};
        \node [style=none] (8) at (3, 3) {};
        \node [style=none] (9) at (5, 2) {};
        \node [style=none] (10) at (3, 2) {};
        \node [style=none] (11) at (20, 2.5) {};
        \node [style=none] (12) at (1, 0) {};
        \node [style=none] (102) at (6, 2.5) {=};
        \node [style=none] (103) at (7, 2.5) {\(0\)};
        \node [style=none] (105) at (13, 3) {};
        \node [style=none] (106) at (9, 5) {};
        \node [style=none] (107) at (11, 3) {};
        \node [style=none] (108) at (9, 2) {};
        \node [style=none] (109) at (11, 2) {};
        \node [style=none] (110) at (13, 0) {};
        \node [style=none] (111) at (8, 2.5) {=};
        \node [style=none] (112) at (15.5, 2.5) {and};
        \node [style=none] (113) at (21, 2.5) {\(=\)};
        \node [style=none] (114) at (22, 2.5) {\(1\)};
      \end{pgfonlayer}
      \begin{pgfonlayer}{edgelayer}
        \draw [very thick, in=-90, out=-90, looseness=1.50] (10.center) to (9.center);
        \draw [very thick, in=270, out=90] (12.center) to (2.center);
        \draw [very thick, in=90, out=-90] (3.center) to (9.center);
        \draw [very thick, in=270, out=90] (10.center) to (8.center);
        \draw [very thick, in=90, out=90, looseness=1.50] (2.center) to (8.center);
        \draw [very thick, in=630, out=270, looseness=1.75] (11.center) to (7.center);
        \draw [very thick, in=90, out=90, looseness=1.75] (7.center) to (11.center);
        \draw [very thick, in=-90, out=-90, looseness=1.50] (109.center) to (108.center);
        \draw [very thick, in=-90, out=90] (110.center) to (105.center);
        \draw [very thick, in=90, out=-90] (106.center) to (108.center);
        \draw [very thick, in=-90, out=90] (109.center) to (107.center);
        \draw [very thick, in=90, out=90, looseness=1.50] (105.center) to (107.center);
      \end{pgfonlayer}
    \end{tikzpicture}
  \end{center}
\end{remark}

\begin{definition}\label{def:sl2-cry}
  We call the Cauchy completion \(\slcrys\) of \(\mathsf{CTL}^{\!\!\op}\) the category of \emph{\(\mathfrak{sl}_{2}(\mathbb{C})\)-crystals}.
  Given \(X, Y \in \slcrys\), we write \(\Hom_{\mathfrak{sl}}(X,Y)\) for the space of morphisms between \(X\) and \(Y\).
\end{definition}
Since \(\mathsf{CTL}\) satisfies Hypothesis~\ref{hypo:eso-enriched}, the tensor product of \(\mathsf{CTL}^{\!\!\op}\) turns \(\slcrys\) into a monoidal category, see \cite[Section~4.1]{alqady-stroiński2025:TemperleyLieb}, raising the question whether internal homs exists and what kind of duality is exhibited by \(\slcrys\).
In order to answer it, we rely on its semisimplicity, defined by a standard generalisation, see \eg \cite{morrison-snyder2012:NonCyclotomicFusionCategories}, of the abelian case discussed in \cite[Definition~1.5.1]{Etingof2015}.

\begin{definition}\label{def:semisimple-category}
  Consider a Cauchy-complete \(\kvect\)-enriched category \(\cat{C}\).
  An object \(x \in \cat{C}\setminus\{0\}\) is \emph{simple} if \(0\) and \(x\) are its only subobjects.

  We call \( \cat{C}\) \emph{semisimple} if every object is a direct sum of simples.
\end{definition}

The next result follows from \cite[Corollary 3.32]{alqady-stroiński2025:TemperleyLieb}.

\begin{proposition}\label{prop:sesi-of-crystals}
  The category \(\slcrys\) is semisimple and there exists a exists a bijection
  \begin{equation*}
    \mathbb{N}_{0}\to \{\text{isomorphism classes of simple objects}\}, \qquad
    n \mapsto [L_{n}]
  \end{equation*}
  satisfying for all \(n,m\in \mathbb{N}_{0}\) the Clebsch--Gordon decomposition identities
  \begin{equation}\label{eq:Clebsch-Gordon}
    L_{0} \cong 1, \qquad L_{1} \cong t, \qquad
    L_{n} \otimes L_{m} \cong \oplus^{\min(n,m)}_{i=0} L_{n+m - 2i}.
  \end{equation}
\end{proposition}
By Lemma~\ref{lem:copowering-presheaves}, \(\widehat{\mathsf{CTL}}\) is copowered over \(\kVect[ \mathbb{C}]\), implying that for any \(X \in \slcrys\) and finite-dimensional \(v\in \kvect[\mathbb{C}]\), we have \(X \boxtimes v \in \slcrys\).
For any \(n\in \mathbb{N}_{0}\), we write
\begin{equation}
  \ev_{n,X} \from L_{n} \boxtimes \Hom_{\mathfrak{sl}}(L_{n}, X) \to X
\end{equation}
for the canonical \emph{evaluation} homomorphism which corresponds under the copower adjunction to \(\id \in \kVect[ \mathbb{C}](\Hom_{\mathfrak{sl}}(L_{n}, X), \Hom_{\mathfrak{sl}}(L_{n}, X))\).
This leads to a canonical decomposition of \(X\) into isotypical components.
\begin{lemma}\label{lem:canonical-decomposition}
  For any object \(X\in \slcrys\) the canonical evaluations induce an isomorphism
  \begin{equation}
    \can \eqdef \sum_{n\in \mathbb{N}_{0}} \ev_{n,X} \from \bigoplus_{n\in \mathbb{N}_{0}} L_{n} \boxtimes \Hom_{\mathfrak{sl}}(L_{n}, X) \to X.
  \end{equation}
\end{lemma}

In order to treat the monoidal unit \(L_{0}\) uniformly with the set of other simples, we define \(L_{-1}\eqdef 0 \in \slcrys\).

\begin{remark}\label{rmk:decompositions}
  Since \(L_{1}\otimes L_{n} \cong L_{n+1} \oplus L_{n-1} \) for all \(n\in \mathbb{N}_{0}\), there are two canonical embeddings \(\iota^{\pm}_{n} \from L_{n\pm 1} \to L_{1}\otimes L_{n}\) whose retracts we denote by \(\pi_{n}^{\pm} \from L_{1}\otimes L_{n} \to L_{n\pm 1}\).
  These allow us to form the morphisms
  \begin{subequations}
    \begin{align}
      \varphi^{+}_{n} \eqdef L_{n} \xrightarrow{\iota^{-}_{n+1}} L_{1} \otimes L_{n+1} \xrightarrow{L_{1} \otimes \iota^{+}_{n}} L_{1}\otimes L_{1} \otimes L_{n}, \\
      \varphi^{-}_{n} \eqdef L_{n} \xrightarrow{\iota^{+}_{n-1}} L_{1}\otimes L_{n-1} \xrightarrow{L_{1}\otimes \iota^{-}_{n}} L_{1} \otimes L_{1} \otimes L_{n}.
    \end{align}
  \end{subequations}
  The intersection\footnote{\,%
    Here, \(\im \varphi^{+}\cap \im \varphi^{-}\) is defined as the pullback of the canonical embeddings of  \(\im \varphi^{+}\) and \(\im \varphi^{-}\) into \(L_{1} \otimes L_{1} \otimes L_{n}\).}
  of the images \(\im \varphi^{+}\cap \im \varphi^{-}\) is the zero object.
\end{remark}

Using the compatibility between the tensor products and copowers established in Lemma~\ref{lem:copowering-presheaves}, we can lift the previous construction to arbitrary objects.

\begin{definition}\label{def:the-maps}
  Given any \(X\in \slcrys \), we define two morphisms
  \begin{equation}\label{eq:the-can-embeddings-tensoring-w-L1}
    \varphi^{\pm}_{X}\from X\to L_{1}\otimes L_{1}\otimes X
  \end{equation}
  via the following commutative diagram.
  \begin{equation*}
    \begin{tikzcd}[ampersand replacement=\&,cramped]
      X \&\&\&\& {L_1 \otimes L_1 \otimes X} \\
      \\
      {\bigoplus_{n \in \mathbb{N}_0} L_n \boxtimes \Hom_{\mathfrak{sl}}(L_n, X)} \&\&\&\& {\bigoplus_{n\in \mathbb{N}_0} L_1 \otimes L_1 \otimes L_n \boxtimes \Hom_{\mathfrak{sl}}(L_n,X)}
      \arrow["{\varphi^{\pm}_X}", from=1-1, to=1-5]
      \arrow["\can^{-1}"', from=1-1, to=3-1]
      \arrow["{\sum_{n\in \mathbb{N}_0}\varphi_n^\pm\boxtimes \Hom_{\mathfrak{sl}}(L_n,X)}", from=3-1, to=3-5]
      \arrow["{L_1 \otimes L_1 \otimes \can}"', from=3-5, to=1-5]
    \end{tikzcd}
  \end{equation*}
\end{definition}

Our next result shows that the above defined maps \(X\to L_{1}\otimes L_{1} \otimes X\) are natural.

\begin{lemma}\label{lem:naturality}
  Consider objects \(X, Y \in \slcrys\).
  For any morphism \(f \from X\to Y\), we have
  \begin{equation*}
    (L_{1}\otimes L_{1} \otimes f)\varphi^{\pm}_{X} = \varphi^{\pm}_{Y}f.
  \end{equation*}
\end{lemma}
\begin{proof}
  To keep our proof concise, we restrict ourselves without loss of generality to the case of \(X\cong L_{n}^{a}\) and \(Y \cong L_{m}^{b}\) for \(n,m,a,b \in \mathbb{N}_{0}\).
  If \(n\neq m\), we have \(f=0\) and the claim follows immediately.
  In case that \(n = m\), consider the following diagram,
  the notation and commutativity of which we establish below:
  \begin{equation*}
    \begin{tikzcd}[ampersand replacement=\&,cramped]
      X \&[-1em]\&\&\&[-1em] Y \\[1em]
      \& {L_n \boxtimes \Hom_{\mathfrak{sl}}(L_n, X)} \&\& {L_n \boxtimes \Hom_{\mathfrak{sl}}(L_n, Y)} \\
      \\
      \&[-10em] {L_1 \otimes L_1 \otimes L_n \boxtimes \Hom_{\mathfrak{sl}}(L_n, X)} \&\& {L_1 \otimes L_1 \otimes L_n \boxtimes \Hom_{\mathfrak{sl}}(L_n, Y) } \\[1em]
      {L_1 \otimes L_1 \otimes X} \&\&\&\& {L_1 \otimes L_1 \otimes Y}
      \arrow["f", from=1-1, to=1-5]
      \arrow["\varphi^{\pm}_{X}"', from=1-1, to=5-1]
      \arrow["\varphi^{\pm}_{Y}", from=1-5, to=5-5]
      \arrow["\ev_{n, X}^{-1}"{description}, from=1-1, to=2-2]
      \arrow["{L_n \boxtimes f'}", dashed, from=2-2, to=2-4]
      \arrow["{\varphi^{\pm}_{n}\boxtimes \Hom_{\mathfrak{sl}}(L_n,X)}"', from=2-2, to=4-2]
      \arrow["\ev_{n, Y}^{-1}"{description}, from=1-5, to=2-4]
      \arrow["{\varphi^{\pm}_{n}\boxtimes \Hom_{\mathfrak{sl}}(L_n,Y)}", from=2-4, to=4-4]
      \arrow["{L_1 \otimes L_1 \otimes L_n \boxtimes f'}"', dashed, from=4-2, to=4-4]
      \arrow["{L_1\otimes L_1 \otimes \ev_{n,X}}"{description}, from=4-2, to=5-1]
      \arrow["{L_1 \otimes L_1 \otimes \ev_{n, Y}}"{description}, from=4-4, to=5-5]
      \arrow["{L_1 \otimes L_1 \otimes f}"', from=5-1, to=5-5]
    \end{tikzcd}
  \end{equation*}
  Due to Lemma~\ref{lem:canonical-decomposition} all diagonal arrows are isomorphisms and it suffices to show that the smaller quadrilaterals commute in order for the outer rectangle to be commutative.

  While the left and right trapezia simply state the definition of the morphisms \(\varphi_{X}^{\pm}\) and \(\varphi^{\pm}_{Y}\), both paths in the central rectangle yield the arrow \(\varphi^{\pm}\boxtimes f'\).
  By Corollary~\ref{cor:copwering-is-faithful}, the map
  \begin{align*}
    \kVect[ \mathbb{C}](\Hom_{\mathfrak{sl}}(L_{n},X), \Hom_{\mathfrak{sl}}(L_{n}, Y))
    & \to \Hom_{\mathfrak{sl}}(L_{n}\boxtimes \Hom_{\mathfrak{sl}}(L_{n},X), L_{n}\boxtimes \Hom_{\mathfrak{sl}}(L_{n}, Y)), \\
    g &\mapsto L_{n}\boxtimes g
  \end{align*}
  is injective and therefore also surjective as the target and source spaces have the same dimension.
  Thus, there exists a unique \(f' \from \Hom_{\mathfrak{sl}}(L_{n},X) \to \Hom_{\mathfrak{sl}}(L_{n},Y)\) such that the top and bottom trapezia commute.
\end{proof}

The homomorphisms \(\varphi^{+}\) and \(\varphi^{-}\) can be assembled into the unit of an adjunction \(\adj{L_{1}\otimes \blank}{L_{1}\otimes \blank }{\slcrys}{\slcrys}\).
\begin{lemma}\label{lem:adjunction}
  There is a natural isomorphism
  \begin{equation}
    \Phi_{X,Y}\from \Hom(L_{1} \otimes X, Y) \to \Hom(X, L_{1}\otimes Y), \qquad g \mapsto (L_{1} \otimes g)(\varphi^{+}_{X} + \varphi^{-}_{X}).
  \end{equation}
\end{lemma}
\begin{proof}
  We fix \(X,Y \in \slcrys\) and consider morphisms \(f\from W \to X\) and \(h \from Y \to Z\).
  A direct computation shows that
  \begin{align*}
    \Phi_{W,Z}(hg(L_{1}\otimes f))
    & = (L_{1} \otimes hg(L_{1}\otimes f))(\varphi^{+}+ \varphi^{-})\\
    & = (L_{1} \otimes h)(L_{1} \otimes g) (L_{1}\otimes L_{1} \otimes f) (\varphi^{+}+ \varphi^{-})\\
    \overset{\ref{lem:naturality}}& =(L_{1} \otimes h)(L_{1} \otimes g) (\varphi^{+}+ \varphi^{-}) f
                                    = (L_{1} \otimes h) \Phi_{X,Y}(g) f.
  \end{align*}
  To show that \(\phi_{X,Y}\) is invertible, we may restrict to \(X=L_{n}\) and \(Y= L_{m}\).

  In case \(n+1=m\), we have that \(\Hom(L_{1}\otimes L_{n}, L_{m}) = \spanset_{ \mathbb{C}}\{\pi_{n}^{+}\}\) and similarly \(\Hom(L_{n}, L_{1}\otimes L_{m}) = \spanset_{ \mathbb{C}}\{\iota_{n+1}^{-}\}\).
  As discussed in Remark~\ref{rmk:decompositions}, \(\pi^{+}_{n}\iota^{-}_{n} = 0\) implying that
  \begin{align*}
    \Phi_{L_n, L_m}(\pi_{n}^{+})
    & = (L_{1} \otimes \pi_{n}^{+})(\varphi^{+}_{n} + \varphi^{-}_{n})
      = (L_{1} \otimes \pi_{n}^{+})(L_{1}\otimes \iota_{n}^{+})\iota_{n+1}^{-}
      + (L_{1} \otimes \pi_{n}^{+})(L_{1}\otimes \iota_{n}^{-})\iota_{n-1}^{+}\\
    & =  (L_{1} \otimes \pi_{n}^{+})(L_{1}\otimes \iota_{n}^{+})\iota_{n+1}^{-}
      = (L_{1} \otimes L_{n+1})\iota_{n+1}^{-}
      = \iota_{n+1}^{-}.
  \end{align*}
  Similarly, one shows that \(\Phi\) is bijective in case \(n-1 = m\).
  Finally, if \(n\pm 1 \neq m\), we have \(\dim \Hom(L_{1}\otimes L_{n}, L_{m})= 0 = \dim \Hom (L_{n}, L_{1} \otimes L_{m})\).
\end{proof}

It was shown in~\cite{etingof-penneys2025:RigidityNonNeglegibleModerateGrowth} that \(\mathfrak{sl}_{2}(\mathbb{C})\)-crystals are a non-rigid \(r\)-category.
We will now prove that they are moreover tensor-representable.

\begin{theorem}\label{thm:sl2-is-trep}
  The category \(\slcrys\) is a non-rigid tensor-representable \(r\)-category.
\end{theorem}
\begin{proof}
  Iterated applications of Lemma~\ref{lem:adjunction} show that for all \(n\in \mathbb{N}_{0}\), there is an isomorphism
  \begin{equation*}
    \Hom(L_{1}^{\otimes n} \otimes X, Y)\cong \Hom(L_{1}^{\otimes n-1} \otimes X, L_{1}\otimes Y )\cong \dots \cong \Hom(X, L_{1}^{\otimes n} \otimes Y)
  \end{equation*}
  natural in \(X\) and \(Y\).

  As \(\mathsf{CTL}^{\!\!\op}\) can be identified with the full subcategory of \(\slcrys\) whose objects are tensor powers of \(L_{1}\), the category \(\mathsf{CTL}^{\!\!\op}\) is right tensor-representable.
  Thus, by Corollary~\ref{cor:pres-ref-dual-struct}, we have that \(\slcrys\) is right tensor-representable.
  By~\cite[Theorem~4.10]{alqady-stroiński2025:TemperleyLieb} we have \( \slcrys\cong \slcrys^{\tensorop}\). Therefore, \(\slcrys\) is also left tensor-representable.

  If \(\slcrys\) were rigid, the discussion after Proposition~\ref{prop:dualisable-conditions-Deligne} shows that \(L_{1}\) would be its own rigid dual.
  By Remark~\ref{rmk:temperley--lieb-as-string-diags} the evaluation and coevaluation morphisms are necessarily given by \(\lambda e \from L_{1} \otimes L_{1} \to L_{0}\) and \(\mu c \from L_{0} \to L_{1}\otimes L_{1}\) for \(\lambda, \mu \in \mathbb{C}\).
  However, since \((L_{1}\otimes e)(c \otimes L_{1})=0\) the snake-identities cannot hold and therefore \(\slcrys\) cannot be rigid.

  We conclude the proof by showing that \(\slcrys\) is an \(r\)-category;
  \ie, that there is equivalence of categories
  \begin{equation*}
    [-,1]_{r} \from \slcrys^{\op} \to \slcrys.
  \end{equation*}
  The above computations show that, on objects, we have \({[X, 1]}_{r} \cong X \otimes 1 \cong X\) implying that \({[-,1]}_{r}\) is essentially surjective.
  It being fully faithful follows from the observation that for any simple object \(L_{n}\), we have \(\dim \Hom(L_{n}, L_{n}) = 1 = \dim \Hom([L_{n},1], [L_{n},1])\).
  Therefore, \(\lambda \id_{L_{n}} \in \Hom(L_{n}, L_{n})\)
  gets sent to \(\lambda \id_{[L_{n},1]} \in \Hom([L_{n}, 1], [L_{n}, 1])\),
  and thus \({[\blank, 1]}_r\) is an equivalence.
\end{proof}

\subsection{Boolean algebras}\label{sec:boolean-algebras}
Lattices abstract the idea of forming unions and intersections, and provide structure to a wide range of algebraic phenomena.
The two binary operations can be used to form a poset-category with the lattice as set of objects.
Assuming distributivity between these operations as well as the existence of maximal and minimal elements, it can be endowed with two monoidal structures.
Boolean algebras additionally require complements to exists.
This translates to a \GV{} structure on the poset-category and transforms one tensor product into the other.

\begin{definition}\label{def:lattice}
  A \emph{lattice} is a set \(L\) with two associative and commutative operations \(\land \from L\times L \to L, \quad (a,b) \mapsto a \land b \qquad \text{and}\qquad
  \lor \from L \times L \to L, \quad (a,b) \mapsto a \lor b\), called \emph{meet} and \emph{join}, which satisfy the \emph{absorption laws}
  \begin{equation}\label{eq: absoption-laws}
    a \lor ( a\land b) = a,\qquad
    a \land (a \lor b) = a, \qquad\quad \text{for all }a, b \in L.
  \end{equation}
\end{definition}

Any lattice \((L, \land, \lor)\) defines a poset via the relation
\begin{equation}\label{eq:poset-from-lattice}
  a \leq b \iff a \lor b = b, \qquad\qquad a, b \in L.
\end{equation}
Conversely, a poset \((P,\leq)\) that admits for any pair of objects \(a,b \in P\) a least upper bound \(a\lor b \in P\) and a greatest lower bound \(a \land b\in P\) is a lattice.

A direct computation shows that an element \(t\in L\) is maximal with respect to the partial order of the lattice \(L\) if and only if \(t \lor a = t\)  for all \(a\in L\).
Analogously, \(s\in L\) being minimal equates to \(s\land a = s\) for any \(a\in L\).
Minimal and maximal elements are unique.
In case they exist, we call \(L\) \emph{bounded}.

\begin{definition}\label{def:Boolean-algebra}
  A \emph{Boolean algebra} is a bounded lattice \((L, \land, \lor)\) satisfying \emph{distributivity}
  \begin{equation}\label{eq:distributivity-equation}
    a \land (b \lor c) = (a \land b) \lor (a \land c) \quad \text{and} \quad
    a \lor (b \land c) = (a \land b) \lor (a \land c) \quad \text{for all }a,b,c\in L,
  \end{equation}
  in which every \(a\in L\) admits a \emph{complement} \(a^{\perp}\in L\) in the sense that
  \begin{equation}\label{eq:lattice-complement-mult}
    a \lor a^{\perp} = t \qquad \text{ and }\qquad
    a \land a^{\perp} = s.
  \end{equation}
\end{definition}

Notice that any of the two distributivity requirements given in Equation~\eqref{eq:distributivity-equation} implies the other.
Furthermore, a direct computation shows that complements are unique.
Therefore, we obtain an involutive map
\begin{equation*}
  {( \blank )}^{\perp} \from L \to L, \qquad\quad a \mapsto a^{\perp}.
\end{equation*}

The maximal element \(t\in L\) of a Boolean algebra \((L, \land, \lor)\) is a unit for \(\land\):
\[
  a \land t = a \land (a \lor a^{\bot}) = a, \qquad \text{for all } a \in L.
\]
An analogous computation shows that its minimal element \(s \in L\) satisfies \(a \lor s = a\).

\begin{example}
  We briefly discuss three examples of Boolean algebras with connections to group and ring theory.
  \smallskip

  \textit{Central idempotents}.
  The set \(C\eqdef \{\,e\in Z(R)\mid e^2=e\,\}\) of central idempotents of a unital ring \(R\) is a Boolean algebra when endowed with the operations
  \begin{gather*}
    \land \from C\times C \to C, \quad e \land f = ef, \qquad
    \lor \from C \times C \to C, \quad e \lor f = e+f -ef, \\
    {( \blank )}^{\perp}\from C \to C, \quad e^{\perp}=1-e.
  \end{gather*}

  \textit{Annihilators in semiprime rings}.
  Consider a commutative semiprime ring \(R\) with \(1\).
  That is, its Jacobson radical is trivial.
  In case \(R\) is furthermore Artinian, this is equivalent to it being semisimple.
  The \emph{annihilator} of an ideal \(I \subset R\) is \(I^{\perp} \eqdef \{\,x\in R \mid xI=0\,\}\);
  it is a radical ideal.
  We define on the set \(\mathrm{Ann}(R)\) of annihilators the maps
  \begin{gather*}
    \land \from \mathrm{Ann}(R) \times \mathrm{Ann}(R)\to \mathrm{Ann}(R), \qquad I \land J = I \cap J, \\
    \lor \from \mathrm{Ann}(R) \times \mathrm{Ann}(R)\to \mathrm{Ann}(R), \qquad I \lor J = {(I + J)}^{\perp}.
  \end{gather*}
  This defines the structure of a Boolean algebra on \(\mathrm{Ann}(R)\).
  In fact, it is shown in~\cite{dube-taherifar2021:LatticeAnnihilatorIdeasSemiprime} that the (right) annihilators of a not necessarily commutative ring form a Boolean algebra if and only if the ring is semiprime.\medskip

  \textit{The subgroup lattice}.
  Let \(H\subseteq G\) be a subgroup of a finite group.
  We write \(H^{\perp}\) for the intersection of all maximal subgroups that do not contain \(H\).
  The minimal group constructed in this manner is the Frattini subgroup \(\Phi(G)=G^{\perp}\) of \(G\).
  Deaconescu, Isaacs, and Wall showed in~\cite{deaconescu-isaacs-walls2011:BooleanLatticeFrattiniSubgroup} that the set
  \begin{equation*}
    \{\, \Phi(G) \subseteq H \subseteq G \mid G= HH^{\perp} \,\},
  \end{equation*}
  partially ordered under inclusion forms a Boolean algebra.
\end{example}

In order to apply the techniques for linear functor categories developed in Section~\ref{sec:tensor-rep-functor-cats} to ordinary functors whose source is the poset-category of a Boolean algebra and whose target is the category of vector spaces, we have to `linearise' ordinary (co)presheaves.

\begin{remark}\label{rmk:linearisation}
  Let \(\cat{C}\) be an ordinary category.
  Its \emph{linearisation} \(\k\cat{C}\) has the same objects and the space of morphisms between two objects \(a, b \in \k\cat{C}\) is \(\k \cat{C}(a, b) \eqdef \spanset_{\k}\cat{C}(a,b)\).
  To define its composition, we note that for any  \(a, b, c \in \k\cat{C}\) there is a unique linear map
  \begin{equation*}
    \blank \circ \blank \from \k\cat{C}(b,c) \otimes_{\k} \k\cat{C}(a,b) \to \k\cat{C}(a,c)
  \end{equation*}
  such that \(g\circ f = gf\)  for all  \(g\in \cat{C}(b,c)\) and  \(f\in \cat{C}(a,b)\).

  Let \(\iota\from \cat{C} \to \k\cat{C}\) be the functor which is the identity on objects and maps any morphism \(f \in \cat{C}(a,b)\) to the corresponding basis vector \(f \in \k\cat{C}(a,b)\).
  For any ordinary functor \(F \from \cat{C} \to \cat{D}\) whose codomain is a \(\k\)-linear category, we obtain a commuting triangle:
  \begin{equation*}
    \begin{tikzcd}[ampersand replacement=\&]
      {\cat{C}} \&\& \cat{D} \\
      \& {\k\cat{C}}
      \arrow["F", from=1-1, to=1-3]
      \arrow["\iota"', from=1-1, to=2-2]
      \arrow["{\exists!G \; (\text{linear})}"', from=2-2, to=1-3]
    \end{tikzcd}
  \end{equation*}
  As a consequence, if we endow the (ordinary) functor category \(\langle \cat{C}, \kVect \rangle\) with the pointwise \(\k\)-linear structure, we obtain an isomorphism of linear categories between \(\langle \cat{C}, \kVect \rangle\) and
  \([\k\cat{C}, \kVect]\).
\end{remark}

Barr showed in \cite{barr79} that the duality of Boolean algebras fits within the framework of \GV{} categories.

\begin{proposition}\label{prop:complement-lattice-form-GV-category}
  The poset-category \(\mathcal{L}\) of a Boolean algebra \((L, \land, \lor)\)  is a left \GV{} category with meet as tensor product, \(t\in L\) as  monoidal unit, and \(s\in L\) as dualising object.
  The dualising functor is induced by the involution \({( \blank )}^{\perp}\from L \to L\).

  If \(L\) is furthermore finite, this leads to a right \GV{} structure on the category \(\langle \cat{L}, \kvect\rangle\) of ordinary functors from \(\cat{L}\) to the category of finite-dimensional \(\k\)-vector spaces with \(\rd*{\k\cat{L}( \blank , s)}\) as dualising object and dualising functor
  \begin{equation*}
    \mathsf{D}_r \from \langle \cat{L}, \kvect\rangle^{\op} \to \langle \cat{L}, \kvect\rangle, \qquad
    F \mapsto {F(\blank^{\perp})}^{*}.
  \end{equation*}
\end{proposition}
\begin{proof}
  A direct consequence of Equation~\eqref{eq:distributivity-equation} is  that for \(a \leq c\) and \(b \leq d\), we have
  \((a \land b) \leq (c \land d)\).
  Thus, since \(\land \from L \times L \to L\) is associative, commutative, and unital, it induces a (symmetric) monoidal structure on \(\cat{L}\).

  To see that \((\cat{L}, s)\) is a  \GV{} category,
  we compute its left internal hom.
  Let us fix a triple of elements \(a,b,c \in L\).
  We have
  \begin{align*}
    b \land a \leq c
    & \iff (b \land a) \lor a^{\bot} \leq c \lor a^{\bot}
      \iff (b \lor a^{\bot}) \land (a \lor a^{\bot}) \leq c \lor a^{\bot} \\
    & \iff b \lor a^{\bot} \leq c \lor a^{\bot}
      \iff b \leq c \lor a^{\perp}
  \end{align*}
  where the last equivalence is due to the fact that \(b = b \lor s \leq b \lor a^{\bot}\)
  and \(b \leq c \lor a^{\perp}\) implying that \(b\lor a^{\perp} \leq (c \lor a^{\perp})\lor a^{\perp}= c \lor a ^{\perp}\).
  It follows that the map \({[\blank, \blank]}_{\ell} \from L \times L \to L\),
  \({[a, b]}_{\ell} \eqdef b \lor a^{\perp}\)
  defines the left internal hom of \(\cat{L}\)
  and the order reversing involution \({[\blank, s]}_{\ell} = {( \blank )}^{\perp} \from L \to L\)
  induces an equivalence of categories \(D_{\ell} \from \cat{L}^{\op}\to\cat{L}\).
  Thus, \((\cat{L}, s)\) is a left \GV{} category in the sense of Remark~\ref{rmk:gv-d-from-D}.

  By Remark~\ref{rmk:linearisation}, we have an equivalence of \(\k\)-linear categories \(\langle \cat{L}, \kvect\rangle\cong [\k\cat{L}, \kvect]\).
  Now, \(\k\cat{L}\) satisfies Hypothesis~\ref{hypo:eso-enriched} and, for \(L\) finite, also
  Hypothesis~\ref{hypothesis:tens-repr-mack}.
  So Proposition~\ref{prop:fun-cat-GV} shows that \(\langle \cat{L}, \kvect\rangle\) has the structure of a right \GV{} category with the specified dualising object and dualising functor.
\end{proof}

Suppose \(L\) is a finite Boolean algebra and \(\cat{L}\) is its poset-category.
As discussed in the first steps of the proof of Lemma~\ref{lem:semisimple-lets-ends-commute}, the category \([\k\cat{L}, \kvect]\) can be identified with the finite-dimensional right modules \(\rmod{A}\) of the path algebra \(A\cong \oplus_{e,f\in L}\k\cat{L}(e,f)\) of the poset associated to \(L\).
The algebras obtained in this manner and their properties are interesting in their own right.

\begin{example}\label{ex:upper-triangular-subalgebra}
  The partial order on the set \(L\eqdef\{\,s,a,b,t\,\}\) displayed in the following diagram defines a Boolean algebra.
  \begin{equation*}
    \begin{tikzcd}[ampersand replacement=\&]
      \& a \\
      s \&\& t \\
      \& b
      \arrow["v", from=1-2, to=2-3]
      \arrow["u", from=2-1, to=1-2]
      \arrow["x"', from=2-1, to=3-2]
      \arrow["y"', from=3-2, to=2-3]
    \end{tikzcd}
  \end{equation*}
  A direct computation shows  that its path algebra \(A\) is isomorphic to a  subalgebra of the \(\k\)-valued \(4\times 4\) upper triangular matrices:
  \begin{equation*}
    A \cong\left\{
      \begin{pmatrix}
        s & u & x & z \\
        0 & a & 0 & v \\
        0 & 0 & b & y \\
        0 & 0 & 0 & t
      \end{pmatrix} \;\middle| \;
      s,u,x,z,a,v,b,y,t \in \k
    \right\}.
  \end{equation*}
\end{example}

Subalgebras of upper triangular matrices were studied by Thrall in~\cite{thrall1948:GeneralisationsQuasiFrobenius} to provide generalisations of \emph{quasi-Frobenius algebras}.
That is, algebras whose projective (right) modules are injective and vice versa.
An algebra \(A\) is called \emph{QF-2} if each indecomposable projective right module and each indecomposable projective left module has a simple socle.

Using the \GV{} structure of Proposition~\ref{prop:complement-lattice-form-GV-category}, we will show that the path algebra of any finite Boolean algebra \(L\) is QF-2.
Hereto we need the following observation.
Taking complements in \(L\) extends to an anti-algebra isomorphism \(\phi \from A \to A\) and the pushforward \({( \blank )}_{\phi} \from \rmod{A} \to \lmod{A}\) is an involutive equivalence of categories.
It maps any left module \(M\) to the right module \(M_{\phi}\) which has the same underlying vector space and is endowed with the action \(a \lact m  \eqdef m \ract \phi(a)\) for all \(m\in M\) and \(a\in A\).

\begin{proposition}\label{prop:qf2-algebras-from-boolean}
  Let \(A\) denote the path algebra of a finite Boolean algebra \(L\).
  The functor
  \begin{equation*}
    \rd*{( \blank_{\phi} )} \from {(\rmod{A})}^{\op} \to \rmod{A}
  \end{equation*}
  is an equivalence of categories;
  in particular \(M\in \rmod{A}\) is projective if and only if \(\rd*{(M_{\phi})}\) is injective.

  Furthermore, \(A\) is QF-2,
  and it being quasi-Frobenius is equivalent to \(|L|=1\).
\end{proposition}
\begin{proof}
  The first statement follows directly from the equivalence of categories between \([\k\cat{L}, \kvect]\) and \(\rmod{A}\) given in the proof of Lemma~\ref{lem:semisimple-lets-ends-commute} and the definition of the dualising functor of \([\k\cat{L}, \kvect]\) stated in
  Proposition~\ref{prop:complement-lattice-form-GV-category}.

  In order to show the second claim, we observe that the finite-dimensional algebra \(A = \oplus_{n \geq 0} A_n\) is graded by the path lengths.
  The elements of \(L\) form a basis of \(A_0\) and correspond to the primitive idempotents.
  That is, non-zero idemponents which cannot be written as a sum of two non-zero orthogonal idempotents.
  Furthermore, the Jacobson radical of \(A\) is \(J(A)= \oplus_{n\geq 1}A_n\).
  The indecomposable projective modules of \(A\) are of the form \(eA\) for \(e\in L\).
  Note that \(eA\) has a vector space basis by paths \([e,y]\) for \(e\leq y \leq t\).
  The socle of \(eA\) is \(\soc(eA)=\{\,m \in eA \mid mJ(A)=0\,\}=\spanset_{\k}\{[e,t]\}\) one-dimensional and therefore simple.
  Using that \(\lmod{A}\cong \rmod{A}\), it follows that \(A\) is a QF-2 algebra.

  For \(|L|=1\), we obtain the trivial (quasi-)Frobenius algebra \(A=\k\).
  Otherwise, there exists an \(e\in L\) such that \(\dim eA\geq 2\).
  A direct computation shows that \(\soc(eA)\cong tA\) and if \(eA\) were injective, the inclusion \(tA\cong\soc(eA)\hookrightarrow eA\) would have a retraction, contradicting the indecomposability of \(eA\).
\end{proof}

\subsection{Mackey functors}\label{sec:mackey-functors}
Let \(G\) be a finite group and write \Sp{G} for the category of isomorphism classes of spans of finite \(G\)-sets.
One can show that each hom-set is a free and finitely-generated commutative monoid, see for example the discussion preceding Lemma 2.1 of~\cite{thevenaz95:mackey}.
In order to obtain a \(\kVect\)-enriched category we  `linearise' \(\Sp{G}\).
That is, we define the category \(\k\Sp{G}\) whose objects coincide with \(\Sp{G}\) and whose hom-spaces are given by
\begin{equation*}
  \k\Sp{G}(x, y) = \k \otimes_{\mathbb{N}} \Sp{G}(x, y), \qquad \qquad \text{ for all } x,y \in \Sp{G}.
\end{equation*}
In other words, any basis of a homomorphism monoid of \(\Sp{G}\) yields a basis of the corresponding vector space in  \(\k\Sp{G}\).

Besides the definition of Mackey functors sketched at the beginning of the section,
there is a more succinct formulation due to Lindner~\cite{lindner76:mackey}.

\begin{definition}\label{def:Mackey-functors}
  Let \(\k\) be a field and \(G\) a finite group.
  The category \(\mky_{\k}(G)\) of \emph{Mackey functors} of \(G\) is given by \([\k\Sp{G},\kVect]\).
\end{definition}

Examples of Mackey functors are plentiful.
Indeed, any representation of \(G\) defines a Mackey functor, see~\cite[Example~53.1]{thevenaz1995:G-algebrasModularRepresentation} as well as~\cite[Propostion~10.1]{panchadcharam07:mackey}.
For an extensive overview, we refer the reader to Chapter 53 of~\cite{thevenaz1995:G-algebrasModularRepresentation}.

One can associate to any finite group \(G\)  a finite-dimensional algebra \(\mathbb{M}G\)—the \emph{Mackey algebra of \(G\)}—whose category of modules is equivalent to \(\mky_{\k}(G)\), see~\cite[Propositions~3.1~and~3.2]{thevenaz95:mackey}.
Its finite-dimensional modules are in correspondence with the (pointwise) finite-dimensional Mackey functors, which we denote by \(\mkyfin_{\k}(G)\).

Given that \(\k\Sp{G}\) is a symmetric monoidal category, \(\mky_{\k}(G)\) is closed symmetric monoidal when equipped with Day convolution as its tensor product.\footnote{\,%
  As such, we will refrain from adding `left' or `right' prefixes to the different duality notions when talking about Mackey functors.%
}
Furthermore, \(\k\Sp{G}\) has a finite dense subcategory   whose objects form a complete set of representatives of transitive \(G\)-sets.
The arguments of Example~\ref{ex:mackey-dense-subcategories} may now be used to show that the tensor product and internal hom of finite-dimensional Mackey functors is finite-dimensional, see~\cite[Section~9]{panchadcharam07:mackey}.
In particular, in this case Hypothesis~\ref{hypothesis:tens-repr-mack} holds.

\begin{proposition}\label{prop:mky-GV-and-rigid}
  Let \(G\) be a finite group and suppose \(\k\) is a field.
  The category \(\mkyfin_\k(G)\) is a \GV{} category with \(\rd*{\mkyfin_\k( \blank,1)}\) as dualising object.
  Its rigid objects are precisely the finitely-generated projective Mackey functors.

  Furthermore, \(\mkyfin_{\k}(G)\) is rigid itself if and only if \(\mathbb{M}G\) is semisimple.
  This is equivalent to \(\characteristic \k\) not dividing the order of \(G\).
\end{proposition}
\begin{proof}
  The fact that \(\mkyfin_{\k}(G)\) is a \GV{} category was proven in~\cite[Theorem~9.2]{panchadcharam07:mackey}.
  Alternatively, we may use that \(\k\Sp{G}\) is a rigid category, see~\cite[Section~2]{panchadcharam07:mackey} and apply
  Proposition~\ref{prop:fun-cat-GV} to obtain that \(\mkyfin_{\k}(G)\)
  is a \GV{} category with \(\rd*{\mkyfin_\k(\blank,1)}\in \mkyfin_{\k}(G)\) as its dualising object.
  Additionally, due to Proposition~\ref{prop:characterisation-of-rigidity}, a Mackey functor \(X \in \mky_{\k}(G)\) is rigid dualisable if and only if it is finitely-generated projective.

  Now let us assume \(\mkyfin_{\k}(G)\) is rigid and recall that it is equivalent to the category \(\lmod{\mathbb{M}G}\) of finite-dimensional modules of the Mackey algebra.
  Since \(\mathbb{M}G\) is finite-dimensional and every object in \(\mkyfin_{\k}G\) is projective by Proposition~\ref{prop:characterisation-of-rigidity}, any submodule of \(\mathbb{M}G\) must have a complement.
  In particular, \(\mathbb{M}G\) is semisimple.
  Conversely, in case \(\mathbb{M}G\) is semisimple, all objects of \(\lmod{\mathbb{M}G}\), and therefore also \(\mkyfin_{\k}(G)\), are projective.
  As they are furthermore finitely-generated, Proposition~\ref{prop:characterisation-of-rigidity} implies the rigidity of \(\mkyfin_{\k}(G)\).

  Corollary~14.4 of~\cite{thevenaz95:mackey}, shows that \(\mathbb{M}G\) is semisimple if \(\characteristic \k\) does not divide \(|G|\).
  On the other hand, the category \(\Rep_{\k}(G)\) of finite-dimensional representations of \(G\) over \(\k\) is a full subcategory of \(\mkyfin_{\k}(G)\), see~\cite[Proposition~10.1]{panchadcharam07:mackey}.
  Thus, if \(\mkyfin_{\k}(G)\) is semisimple, so is \(\Rep_{\k}(G)\),
  implying that \(\characteristic \k\) does not divide \(|G|\)
  by Maschke's theorem.
\end{proof}
Rigidity for Mackey functors was also discussed by Bouc, see~\cite[Lemma~2.2]{bouc05:mackey}.
The previous lemma gives an alternative proof for the sketched argument that the classes of rigid and finitely-generated projective Mackey functors coincide.

\subsection{Representations and 2-groups}\label{sec:2-groups}
Another kind of example for \GV{} structures on abelian \(\k\)-linear functor categories arises from studying (strict) 2-groups, which can be identified with crossed modules.
Our recollection of their definition and basic properties is based on~\cite{wagemann2021:CrossedModules}.

\begin{definition}\label{def:str-2-grp}
  A \emph{strict \(2\)-group} is a (small) groupoid \(\cat{G}\) endowed with a strict monoidal structure such that the monoid of objects \((\Ob(\cat{G}), \otimes,1)\) is a group.
\end{definition}

Every strict 2-group is a rigid category;
the left and right dual of an object \(g \in \cat{G}\) is given by its inverse \(g^{-1} \in \cat{G}\).

Let \(G\) be the group of objects of \(\cat{G}\) and \(H=\cat{G}(1, \blank)\) the set of arrows which start at the monoidal unit \(1\in G\).
Writing \(t \from H \to G \) for the \emph{target map}, we define a group structure on \(H\) via the multiplication
\begin{equation*}
  h' h \eqdef (h' \otimes \id_{t(h)})h, \qquad \text{for all }h,h' \in H.
\end{equation*}
For any \(g\in G\) and \(h\in H\) there is a unique \(h' \in H\) such that \(\id_g\otimes h= h'\otimes \id_g\).
This induces a map \(\alpha \from G \to \Aut(H)\) and turns the quadruple \((G, H, t, \alpha)\) into a crossed module as defined below.

\begin{definition}\label{def:crossed-module}
  A \emph{crossed module} is a quadruple \((G,H,t,\alpha)\) comprising two groups \(G,H\) and two group homomorphisms \(t \from H \to G\) and \(\alpha \from G \to \Aut(H)\), satisfying
  \begin{equation} \label{eq:cross-mod-def-rel}
    t(\alpha(g) h) = gt(h)g^{-1}, \qquad \alpha(t(l))h = lhl^{-1}, \qquad\qquad
    \text{for all \(g \in G\) and \(h,l \in H\)}.
  \end{equation}
\end{definition}

Conversely, any crossed module gives rise to a strict 2-group.

\begin{remark}\label{rmk:crossed-modules-and-2-groups}
  Any crossed module \((G,H,t,\alpha)\) defines a strict \(2\)-group \(\cat{G}\):
  \begin{equation}\label{eq:2-grp-from-crossed-module}
    \begin{gathered}
      \Ob\cat{G}=G ,\qquad\qquad \cat{G}(g,g') =\{(h,g) \in H \times G \mid t(h)g=g' \}, \\
      (h', t(h)g)\circ(h,g)= (h'h,g), \qquad \qquad (h,g)\otimes (h', g')= (h \; \alpha(g)h', g g').
    \end{gathered}
  \end{equation}
  The left and right dualising functors coincide and map any object \(b \in G\) to \(\ld{b}=b^{-1}=\rd{b}\).
  Given a morphism \(f=(h,b)\from b \to c\), we have
  \begin{equation} \label{eq:2-grps-dualising-fun-on-morphs.}
    \ld{f}  =(\alpha(c^{-1})h, c^{-1}) = (\alpha(b^{-1})h, c^{-1}) = \rd{f}.
  \end{equation}
\end{remark}

Indeed, there is an equivalence of categories between crossed modules and strict 2-groups as is shown for example in~\cite[Theorem~1.9.2]{wagemann2021:CrossedModules}.

\begin{example}\label{ex:2-groups-in-algebra-and-rep-theo}
  There are ample examples of crossed modules throughout the literature.
  The following short discussion highlights their versatility and emphasises their applications in various  areas of mathematics.
  \smallskip

  \textit{Hopf--Galois extensions and skew braces}.
  The \emph{holomorph} of a  group \(H\) is the semidirect product \(H \rtimes \Aut(H)\).
  The two group homomorphisms \(t \from H \to \Aut(H)\), \(t(l)h = lhl^{-1}\)and \(\alpha=\id\from \Aut(H) \to \Aut(H)\) turn \((\Aut(H), H, t,\id_{\Aut(H)})\) into a crossed module.
  As explained in the Introduction of~\cite{byott2024:insolubleHolomorph}, Section~2 of~\cite{childs-greither-keating-et-al2021:HopfGalois}, and~\cite[Section~2]{bachiller2016:CounterexampleConjectureBraces}, holomorphs can be used to classify certain Hopf--Galois extensions as well as skew-braces.
  \smallskip

  \textit{Representation theory of finite groups}.
  A  short exact sequence of groups
  \begin{equation*}
    \begin{tikzcd}[ampersand replacement=\&]
      0 \& C \& E \& G \& 0
      \arrow[from=1-1, to=1-2]
      \arrow["\iota", from=1-2, to=1-3]
      \arrow["t", from=1-3, to=1-4]
      \arrow[from=1-4, to=1-5]
    \end{tikzcd}
  \end{equation*}
  such that \(C\) embeds into the centre of \(E\) is called a \emph{central extension} of \(G\).
  Let \(\kappa \from G\to E\) be a set-theoretical section of \(t\from E\to G\).
  A direct computation shows that the map \(\alpha\from G \to \Aut(E)\), \(\alpha(g)e=\kappa(g)e{\kappa(g)}^{-1}\) is independent of the choice of the section, a homomorphism of groups, and the quadruple \((E,G,t,\alpha)\) constitutes a crossed module, see~\cite[Section 1.3]{wagemann2021:CrossedModules}.

  As discussed in~\cite[Chapter~1]{hoffman-humphreys1992:ProjectiveRepresentationsSymmetricGroup}, any projective representation \(\rho \from G \to \PGL_{\mathbb{C}}(V)\) of a finite group \(G\) can be lifted to a linear representation \(\varrho \from E \to \GL(V)\) of a certain central extension \(E\) of \(G\).
  Projective representations themselves are studied for example in the context of representation theory of semidirect products, see~\cite{ceccherini-silberstein-scarabotti-tolli2022:RepresentationFiniteGroupExtensions}.
\end{example}

Let \((G,H,t,\alpha)\) be a crossed module whose associated \(2\)-group we denote by \(\cat{G}\).
Using the results of Section~\ref{sec:tensor-rep-functor-cats}, we determine a \GV{} structure on the category  \(\langle\cat{G}, \kvect\rangle\) of (ordinary) functors between \(\cat{G}\) and \(\kvect\).

Hereto, we need to  analyse the structure of \(\cat{G}\) in more detail.
Due to Equation~\eqref{eq:cross-mod-def-rel}, the image \(K \eqdef\im t \subset G\) is a normal subgroup of \(G\) and the kernel \(L \eqdef \ker t\) is normal and central in \(H\).
The latter follows from the identity
\begin{equation*}
  lhl^{-1} =\alpha(t(l))h = \alpha(e)h = h \qquad\quad \text{ for all } l\in L, h\in H.
\end{equation*}
For any \(k=t(h) \in K\) and \(l\in L\), we get \(\alpha(k)l = \alpha(t(h))l = hlh^{-1}= l\) and therefore, there is a unique homomorphism of groups
\begin{equation}\label{eq:restriction-of-action-by-automorphisms}
  \overline{\alpha} \from Q \eqdef G/K \to \Aut(L), \qquad \text{with } \overline{\alpha}([g])l=\alpha(g)l \in L \subset H \text{ for all } g\in G, l\in L.
\end{equation}

The connected components of \(\cat{G}\) are in bijection with the elements of \(Q\).
Given any element \(g\in G\), we write \( \cat{G}_{g}\) for the maximal full connected pointed subgroupoid of \(\cat{G}\) whose distinguished object is \(g\).
The set of objects of \(\cat{G}_{g}\) is \(Kg\) and its morphisms correspond to \(H \times Kg\).

\begin{lemma}\label{lem:canonical-inclusion}
  Let \(g\in G\) and write \(\mathbf{B}(L)\) for the category with a single object \(g\) whose endomorphisms are given by \(L\).
  The canonical inclusion \(\Theta_g \from \mathbf{B}(L) \to \mathcal{G}_g\) sending \(g\) to \(g\) and \(l \from g\to g\) to \((l,g) \from g\to g\) is essentially surjective, fully faithful, and thereby an equivalence of categories.
\end{lemma}
To specify a quasi-inverse, we fix a set-theoretical section \(\iota \from K \to H\) of  the surjective homomorphisms of groups \(t\from H \to K\) and define \(\Psi_g \from \cat{G}_g\to \mathbf{B}(L)\) by mapping each object to \(g\) and any morphism \((h,kg)\) to \(\iota({(t(h)k)}^{-1}) h \iota(k)\).
Using that \(\langle \mathbf{B}(L), \kvect\rangle\) can be identified with the category \(\Rep(L)\) of finite-dimensional representations of \(L\), the pushforwards of \(\Theta_g\) and \(\Psi_g\) establish an equivalence of categories
\begin{equation}\label{eq:pushforward-equivalence-of-cats}
  \begin{aligned}
    \Theta_{g}^{*} \colon \langle \cat{G}_g, \kvect\rangle &\rightleftarrows \Rep(L) \cocolon \Psi_g^{*} \\
    F &\mapsto (Fg,\rho_{F}) \qquad \text{where } \rho_F(l) = F(l,g) \\
    \left\{\begin{aligned}
      kg &\mapsto M\\
      (hk,g) &\mapsto \rho\left(\iota({(t(h)k)}^{-1}) h \iota(k)\right)
    \end{aligned}\right\} & \longmapsfrom (M, \rho)
  \end{aligned}
\end{equation}

We now define
\begin{equation} \label{eq:rep-cat-assoc-to-2-grp}
  \Rep_{\cat{G}}(L) \eqdef \bigoplus_{q\in Q}{\Rep(L)}_{q}, \qquad \qquad \text{where } {\Rep(L)}_{q} = \Rep(L).
\end{equation}
Given \(M\in \Rep_{\cat{G}}(L)\), we write \(M_q\) for its homogeneous component of degree \(q\in Q\).

Note that in case \(K\) has finite index in \(G\), we have \(\langle \cat{G}, \kvect\rangle \cong \bigoplus_{[g]\in Q}\langle \cat{G}_g, \kvect \rangle\).

\begin{lemma}\label{lem:structure-of-G-funs}
  Let \(\cat{G}\) be a strict \(2\)-group corresponding to a crossed module \((G, H, t, \alpha)\) and assume that \(K=\im t\subset G\) has finite-index.
  Then \(\Rep_{\cat{G}}(L) \simeq \langle \cat{G}, \kvect\rangle\).
\end{lemma}

The tools developed in Section~\ref{sec:tensor-rep-functor-cats} allow us to endow \(\langle\cat{G}, \kvect\rangle\) with the structure of a \GV{} category which we will transfer to \( \Rep_{\cat{G}}(L)\).
The tensor product obtained in this manner will permute the homogeneous components and `twist' the action of \(L\) by virtue of \(\overline{\alpha}\).
In this regard, we introduce the following notation:
for any \(q\in Q\) and \(M \in \Rep(L)\), we denote by \(M^{\overline{\alpha}(q)}\) the representation of \(L\) whose underlying vector space is \(M\), endowed
with the action \(l \blact m = \overline{\alpha}(q)l \lact m\) for all \(l \in L\) and \(m \in M\).

\begin{proposition}\label{prop:2-groups-tensor-product}
  Suppose \((G,H,t,\alpha)\) is a crossed module with \(G\) and \(H\) finite and let \(\cat{G}\) be its associated strict \(2\)-group.
  The category \(\Rep_{\cat{G}}(L)\) is a right \(r\)-category.
  Its tensor product and dualising functor are defined by the assignments
  \begin{equation} \label{eq:gv-structure-on-reps}
    M\otimes N  = (M\otimes_{\k L} N^{\overline{\alpha}(p^{-1})})\in {\Rep_{\cat{G}}(L)}_{pq}, \qquad
    D(M) = \rd*{(M^{\overline{\alpha}(p)})} \in {\Rep_{\cat{G}}(L)}_{p^{-1}}
  \end{equation}
  for \(p, q \in Q\) and \(M \in {\Rep_{\cat{G}}(L)}_{p}\), \(N \in {\Rep_{\cat{G}}(L)}_{q}\).
\end{proposition}
\begin{proof}
  As \(G\) and \(H\) are finite, Hypotheses~\ref{hypo:eso-enriched} and \ref{hypothesis:tens-repr-mack} hold for \(\langle\cat{G},\kvect\rangle \cong [\k\cat{G}, \kvect]\).
  By Proposition~\ref{prop:fun-cat-GV}, \([\k\cat{G}, \kvect]\) can be endowed with the structure of a right \GV{} category with \(\rd*{\k\cat{G}( \blank , 1)}\) as dualising object;
  write \(R \from \langle\cat{G},\kvect\rangle^{\op} \to \langle\cat{G},\kvect\rangle\) for its dualising functor.

  By Lemma~\ref{lem:structure-of-G-funs}, there are \(\k\)-linear equivalences \(\adj{\widehat{\Theta}}{\widehat{\Psi}}{[\k\cat{G}, \kvect]}{\Rep_{\cat{G}}(L)}\) which are determined for each homogeneous component \([\k\cat{G}_{[g]}, \kvect]\) by the  pushforwards of the equivalences \(\adj{\Theta_g}{\Psi_{g}}{\langle\cat{G}_g, \kvect\rangle}{\mathbf{B}(L)}\).
  Using the formulas of Equation~\eqref{eq:pushforward-equivalence-of-cats}, we can therefore explicitly transfer the \GV{} structure of \([\k\cat{G}, \kvect]\) to \(\Rep_{\cat{G}}(L)\).

  The dualising object is mapped to \(\Theta(\rd*{\k\cat{G}( \blank, 1)})=\rd*{\k L} \cong \k L\in {\Rep_{\cat{G}}(L)}_{[1]}\), where we used that  \(\k L\) is a Frobenius algebra for the last equality.
  The tensor product and right dualising functor of \(\Rep_{\cat{G}}(L)\) are given by
  \begin{gather*}
    \otimes \eqdef \widehat{\Theta}\ostar(\widehat \Psi \times \widehat \Psi) \from \Rep_{\cat{G}}(L) \times \Rep_{\cat{G}}(L) \to \Rep_{\cat{G}}(L), \\
    D \eqdef \widehat{\Theta}\mathsf{R}\widehat \Psi^{\op} \from {\Rep_{\cat{G}}(L)}^{\op} \to \Rep_{\cat{G}}(L).
  \end{gather*}
  In order to compute it explicitly, we consider two elements \(p,q\in Q\) as well as two representations \((M, \rho)\in {\Rep_{\cat{G}}(L)}_{p}\) and \((N, \tau) \in {\Rep_{\cat{G}}(L)}_{q}\).
  The Day convolution of the functors corresponding to \(M\) and \(N\) is computed via a certain colimit and \(\widehat{\Theta}\), as an equivalence of categories, commutes with this colimit.
  Thus, by Equation~\eqref{eq:day-convolution-simp-right}, and the definition of coends, the homogeneous component \({(M\otimes N)}_x\) of degree \(x\in Q\) is the coequaliser of the diagram
  \begin{equation*}
    \begin{tikzcd}
      \displaystyle{\coprod_{r, s \in Q}}  \k L_{r,s} \otimes_{\k} M_r \otimes_{\k} N_{s^{-1}x} & \displaystyle{\coprod_{r\in Q}  M_r \otimes_{\k} N_{r^{-1}x},}
      \arrow[shift right=1, from=1-1, to=1-2]
      \arrow[shift left=1, from=1-1, to=1-2]
    \end{tikzcd}
  \end{equation*}
  where \(L_{r,s}=L\) if \(r=s\) and the empty set otherwise.
  Its parallel morphisms are
  \begin{equation*}
    \begin{aligned}
      l\otimes_{\k} m\otimes_{\k} n &\mapsto  \rho(l)m \otimes_{\k} n, \\
      l\otimes_{\k} m\otimes_{\k} n &\mapsto  m \otimes \tau(\rd{l})n = m \otimes \tau(\overline{\alpha}(s^{-1})l)n,
    \end{aligned}
  \end{equation*}
  where \(l \in L_{r,s}\), \(m \in M_{r}\) and \(n \in N_{s^{-1}x}\).

  In  order to compute \(D(M)\), we note that \(\mathsf{R}(F)= \rd*{F(\ld{\blank})}\) for any \(F \in [\k\cat{G}, \kvect]\).
  Thus, given \(x\in Q\), we compute
  \begin{equation*}
    {D(M)}_x =
    \begin{cases}
      \rd*{M} & x^{-1} = p, \\
      \{0\} & \text{otherwise}.
    \end{cases}
  \end{equation*}
  Accordingly, the action of any \(l \in L\) on \(\rd*{M}\) is given by \(\rd*{(\rho(\overline{\alpha}(p)l))}\from \rd*{M}_{p} \to \rd*{M}_{p}\).
\end{proof}

Finally, let us state necessary and sufficient conditions for \(\Rep_{\cat{G}}(L)\) to be rigid.
\begin{proposition}\label{prop:2-grp-GV-and-rigid}
  Let \((G,H,\alpha,t)\) be a crossed module with \(G\) and \(H\) finite and write \( \cat{G}\) for its associated strict 2-group.
  The category \(\Rep_{\cat{G}}(L)\) is rigid if and only if \(\characteristic \k\) does not divide the order of \(L= \ker t\).
\end{proposition}
\begin{proof}
  It follows from Proposition~\ref{prop:characterisation-of-rigidity} that \(\Rep_{\cat{G}}(L)\) is rigid if and only if all of its objects are finitely-generated and projective.
  As \(\Rep_{\cat{G}}(L)\) is a direct sum of finitely many copies of \(\Rep(L)\) and every object of \(\Rep(L)\) is finitely-generated, this is the case if and only if all objects of \(\Rep(L)\) are projective.
  The latter is, in turn, equivalent to \(\k L\) being semisimple which, by Maschke's theorem, corresponds to \(\characteristic \k\) and \(|L|\) being coprime.
\end{proof}
\def\bibfont{\footnotesize}
\bibliographystyle{alpha}
\bibliography{main}

\end{document}